\documentclass[11pt, reqno]{amsart}

\usepackage{amssymb}
\usepackage{hyperref} 
\usepackage{mathrsfs,bbold}
\usepackage{amsmath}
\usepackage{amsthm}
\usepackage{enumerate}
\usepackage{dsfont}
\usepackage{color}
\usepackage[a4paper]{geometry}
\usepackage[utf8]{inputenc}
\usepackage{stmaryrd}
\usepackage{titlesec}
\usepackage{mathabx}
\usepackage{appendix}
\usepackage{todonotes, fancyhdr}
\usepackage{chngcntr}
\usepackage{cancel}

\usepackage{setspace}
\renewcommand{\baselinestretch}{0.99}

\usepackage[all]{xy}
\numberwithin{subsection}{section}
\numberwithin{subsubsection}{subsection}
\numberwithin{equation}{section} 

\newenvironment{Dem}[1][\unskip]{%
    \begin{list}{\hspace{0.5cm}{\sf \textbf{Proof #1 --}}}{%
        \setlength{\topsep}{0pt}%
        \setlength{\leftmargin}{0pt}%
        \setlength{\rightmargin}{0pt}%
        \setlength{\listparindent}{0pt}%
        \setlength{\itemindent}{0pt}%
        \setlength{\parsep}{0pt}%
        \addtolength{\leftmargin}{20pt}%
        \addtolength{\rightmargin}{0pt}%
    } \item }{\hfill $\rhd$\end{list}\smallskip}

\renewcommand\thesection       {\arabic{section}}
\renewcommand\thesubsection    {\thesection{\boldmath $.$}\arabic{subsection}}
\renewcommand\thesubsubsection    {\thesection{\boldmath $.$}\arabic{subsection}{\boldmath $.$}\arabic{subsubsection}}

\titleformat{\section}[block]
{\filcenter\normalfont\sffamily\bfseries\Large}
{{\hspace{-0.7cm}}\thesection \hspace{0.2em} --\vspace{0.3cm}}{0.5em}{}

\titleformat{\subsection}[block]
{\filcenter\normalfont\sffamily\bfseries\large}  						  
{\hspace{-0.7cm}\thesubsection \hspace{0.5em} \vspace{0.3cm}}{.5em}{}  
\titlespacing{\subsection}{-0pc}{1.5ex plus .1ex minus .2ex}{0pc}

\titleformat{\subsubsection}[block]
{\normalfont\sffamily\bfseries}					  
{\thesubsubsection \vspace{0.3cm}}{.5em}{}  
\titlespacing{\subsection}{-0pc}{1.5ex plus .1ex minus .2ex}{0pc}

\newtheoremstyle{mystyle}
{3pt}               
{3pt}               
{\it }                      
{}                      
{\sffamily\bfseries}             
{}                      
{0.5em}                 
{#1 #2{\Large$.$}  }

\theoremstyle{mystyle}

\newtheorem{thm}{Theorem}
\newtheorem*{thm*}{Theorem}

\newtheorem{prop}[thm]{\hspace{-0.17cm} {Proposition}}
\newtheorem{defn}[thm]{ \hspace{-0.33cm} {Definition}}

\newtheorem*{rem*}{\hspace{-0.15cm} {Remark}}

\newtheoremstyle{mystyle2}
{3pt}               
{3pt}               
{\it }                      
{}                      
{\sffamily\bfseries}             
{}                      
{0.5em}                 
{\llap{#2 }#1{\hspace{0.2cm}--}}

\theoremstyle{mystyle2}

\newtheorem*{definition*}{Definition}
\newtheorem*{theorem*}{Theorem}
\newtheorem*{Remark*}{Remark}
\newtheorem*{lem*} {Lemma}
\newtheorem*{defn*} {Definition}
\newtheorem*{prop*} {Proposition}
\newtheorem*{cor*} {Corollary}

\newcommand{\IR}{\mathbb{R}}

\newcommand{\IT}{\mathbb{T}}

\newcommand{\dst}{\displaystyle}
\newcommand{\drm}{\mathrm d}

\newcommand{\CD}{\mathcal D}

\newcommand{\CP}{\mathcal P}

\newcommand{\CM}{\mathcal M}
\newcommand{\CQ}{\mathcal Q}
\newcommand{\mcC}{\mathcal C}

\newcommand{\CG}{\mathcal G}

\newcommand{\SA}{\mathscr{A}}

\newcommand{\SL}{\mathscr{L}}

\renewcommand{\P}{\mathsf{P}}
\newcommand{\PT}{\widetilde{\P}}
\newcommand{\PI}{\mathsf{\Pi}}

\newcommand{\DD}{\mathsf{D}}
\newcommand{\DC}{\mathsf{C}}
\newcommand{\DR}{\mathsf{R}}

\newcommand{\DL}{\mathsf{L}}

\newcommand{\DV}{\mathsf{V}}

\newcommand{\Mi}{\tau}
\newcommand{\Mj}{\sigma}
\newcommand{\Mk}{\gamma}

\newcommand{\ssk}{\smallskip}

\renewcommand{\epsilon}{\varepsilon}

\newcommand\bbR{\mathbb{R}}

\newcommand{\CC}{\mathcal{C}}

\newcommand{\mcM}{\mathcal{M}}

\newcommand{\mcQ}{\mathcal{Q}}
\newcommand\mcS{\mathcal{S}}

\begin{document}

\vspace*{3ex minus 1ex}
\begin{center}
{\Huge\sffamily{Paracontrolled calculus for quasilinear singular PDEs  \vspace{0.5cm}}}
\end{center}

\begin{center}
{\sf I. BAILLEUL\footnote{I.Bailleul thanks the U.B.O. for their hospitality, part of this work was written there. 


} and A. MOUZARD}
\end{center}

\vspace{1cm}

\begin{center}
\begin{minipage}{0.8\textwidth}
\renewcommand\baselinestretch{0.7} \scriptsize \textbf{\textsf{\noindent Abstract.}} We develop further in this work the high order paracontrolled calculus setting to deal with the analytic part of the study of quasilinear singular PDEs. A number of continuity results for some operators are proved for that purpose. Unlike the regularity structures approach of the subject by Gerencser \& Hairer and Otto, Sauer, Smith \& Weber, or Furlan and Gubinelli' study of the two dimensional quasilinear parabolic Anderson model equation, we do not use parametrised families of models or paraproducts to set the scene. We use instead infinite dimensional paracontrolled structures that we introduce here.
\end{minipage}
\end{center}

\vspace{1cm}

\section{Introduction}
\label{SectionIntroduction}

This work is dedicated to the study of the quasilinear singular partial differential equation (PDE)
\begin{equation}   \label{EqQgPAM}
\partial_tu-d(u) Au = f(u) \zeta,
\end{equation}
where $\zeta$ stands for a spacetime noise of parabolic H\"older regularity $\alpha-2$, with $2/5<\alpha<1/2$, with a real-valued unknown $u$ defined on a $3$-dimensional closed Riemannian manifold $M$, and $A$ is an elliptic operator on $M$ in H\"ormander form 
$$
A = \sum_{i=1}^\ell A_i^2,
$$ 
for smooth vector fields $A_i$ on $M$, and $\ell\geq 3$. The function $d$ -- for diffusivity, is supposed to be smooth enough and to take its values in a compact set of $(0,+\infty)$. We assume here for simplicity that the initial condition $u_0$ in equation \eqref{EqQgPAM} is regular enough to treat the free propagation of the initial condition as a remainder term and avoid the technical use of weighted norms. The reader acquainted with the results of Bailleul and Bernicot's work \cite{BB2} on the high order paracontrolled calculus will see that our method for the study of equation \eqref{EqQgPAM}, and the tools introduced along the way, give a direct access to the analysis of the quasilinear generalised (KPZ) equation
$$
\partial_tu-d(u)\partial_x^2u = f(u) \zeta + g(u)\vert \partial_xu\vert^2,
$$
or any other quasilinear version of parabolic semilinear equations, or systems of equations, that can be studied within the setting of the high order paracontrolled calculus.

\ssk

Paracontrolled calculus was introduced in Gubinelli, Imkeller and Perkowski' seminal work \cite{GIP} as a first order `expansion machinery' for the study of a number of singular PDEs. Despite the first order limitation, the paracontrolled approach to the study of singular PDEs has been very successful, as Gubinelli and Perkowski's works \cite{KPZReloaded, GPGenerator} on the KPZ and stochastic Burgers equations, Catellier-Chouk, Mourrat-Weber and Gubinelli \& co-authors works \cite{CatellierChouk, MourratWeber1, MourratWeber2, BarashkovGubinelli, HofmanovaGubinelli} on the $\Phi^4$ scalar equation from quantum field theory, the works \cite{AllezChouk, ChoukvZuijlem} of Chouk and co-authors on the spectral theory for the $2$-dimensional Laplacian with white noise potential, and the very recents works on hyperbolic singular PDEs \cite{GKO, ORT}, testify, amongst others. The scope of the first order paracontrolled calculus was much extended in \cite{BB1, BB2, BB3}, and the high order paracontrolled calculus offers now a convenient setting for the study of a whole class of singular parabolic PDEs, in diverse geometric settings.

\ssk

The study of quasilinear singular PDEs was launched by the works \cite{OttoWeber} of Otto and Weber, \cite{FurlanGubinelli} of Furlan and Gubinelli, and \cite{BDH} of Bailleul, Debussche and Hofmanov\'a, that all appeared within a few months. Interestingly, each of these works used different methods to tackle the same equation: The $2$-dimensional quasilinear parabolic Anderson model equation. Otto and Weber introduced a rough paths flavoured variant of regularity structures, Furlan and Gubinelli introduced a variant of the first order paracontrolled calculus using paracomposition operators instead of paraproducts, while Bailleul, Debussche and Hofmanov\'a showed that the original first order paracontrolled calculus is sufficient to prove well-posedness of the equation on a small time interval. Gerencs\'er and Hairer then showed in \cite{GerencserHairer} that the study of a whole class of quasilinear singular parabolic PDEs can be done in the setting of regularity structures, in the above regime for the regularity exponent $\alpha$, giving results way beyond the scope of what was proved in \cite{OttoWeber, FurlanGubinelli, BDH} and Otto, Sauer, Smith and Weber's followup work \cite{OSSW}. The only caveat to their remarkable results is the fact that their formulation of the quasilinear equation does not allow for a clean treatment of the renormalisation problem; it is thus unclear at the moment that renormalised equations are `local', like in the semilinear setting \cite{BCCH18, BHZ}. See however Gerencs\'er's recent work \cite{Gerencser} for a first result in this direction. 

\ssk

By adding a few results to the toolkit of the high order paracontrolled calculus \cite{BB3}, we are able to prove a local in time well-posedness result for equation \eqref{EqQgPAM}, with the same line of attack as in \cite{BDH}. The method works mutatis mutandis for the study of the quasilinear generalised (KPZ) equation or the quasilinear version of the geometric stochastic heat equation. The present work is purely analytical and does not consider the problem of renormalisation. This amounts here to assuming that a sequence of multilinear functions of the noise are given a priori as elements of their natural spaces, with natural bounds on their norms.

\ssk

$\bullet$ Let $\overline{u_0}$ be a smooth function close enough to $u_0$ in $C^{4\alpha}$ -- this will be quantified later, in Theorem \ref{ThmMain}. Our starting point consists in defining the solution-independent operator
$$
L := -\sum_{i=1}^\ell V_i^2, \qquad V_i := \sqrt{d(\overline{u_0})}\,A_i,
$$
and rewriting equation \eqref{EqQgPAM} under the form of an evolution equation
\begin{equation} \begin{split}   \label{EqReformulation}
\mathscr{L}u :=& \; (\partial_t+L)u   \\
=&\; f(u)\zeta + \big(d(u)-d(u_0)\big)Au + \sum_{i=1}^\ell A_i\big(d(\overline{u_0})^{1/2}\big)V_iu   \\
=& \; f(u)\zeta - d(u_0)^{-1}\big(d(u)-d(\overline{u_0})\big)Lu + \sum_{i=1}^\ell \Big(1- d(\overline{u_0})^{-1}\big(d(u)-d(\overline{u_0})\big)\Big) A_i\big(d(\overline{u_0})^{1/2}\big)V_iu   \\
=&\hspace{-0.12cm}: f(u)\zeta + \varepsilon(u,\cdot)Lu + \sum_{i=1}^\ell a_i(u,\cdot)V_iu,
\end{split}\end{equation}
involving the solution-independent operator $L$. The nonlinear term 
$$
\varepsilon(u,\cdot)Lu = d(\overline{u_0})^{-1}\big(d(u)-d(\overline{u_0})\big)Lu
$$
in the right hand side still involves a second order term, a feature of quasilinear equations. (The dot sign in $\varepsilon(u,\cdot)$ stands for the dependence on $x\in M$ of $\varepsilon$, via $d(\overline{u_0})$.) In the spirit of the study \cite{BDH} of the $2$-dimensional semilinear parabolic Anderson model equation, we are able to define a paracontrolled structure and formulate the equation as a fixed point for a contracting map defined on this structure, using the fact that $u$ stays close to $u_0$ on a small time interval. This is our main result, stated in Theorem \ref{ThmMain}. As a guide for the reader, we describe now the twist on the paracontrolled setting that we use to handle quasilinear equations. 

\ssk

Following \cite{BB2}, one can associate to the differential operator $L$ a paraproduct $\sf P$, and its companion paraproduct $\widetilde{\sf P}$, intertwined to $\sf P$ by the relation
\begin{equation}   \label{EqIntertwining}
\mathscr{L}^{-1}\circ{\sf P} = \widetilde{\sf P}\circ \mathscr{L}^{-1}.
\end{equation}
A resonant operator $\sf \Pi$ is also constructed from $L$. The basic mechanics of the paracontrolled approach to semilinear singular PDEs 
is best illustrated on the model case of the $2$-dimensional parabolic Anderson model equation 
\begin{equation}  \label{EqPAM}
\mathscr{L}u = u\zeta = {\sf P}_u\zeta + {\sf P}_\zeta u + {\sf \Pi}(u,\zeta),
\end{equation}
with null initial condition. One has almost surely $\zeta$ in the parabolic H\"older space $\mcC^{\alpha-2}$, for $\alpha$ any positive real number strictly smaller than $1$. Whereas the above paraproduct terms always make sense for arguments in H\"older spaces of positive or negative exponents, the resonant term is well-defined only if the sum of the H\"older regularity exponents of $u$ and $\zeta$ add up to a positive real number. With $\zeta$ of H\"older regularity $\alpha-2$ and $\alpha<1$, one has $\alpha+(\alpha-2)<0$, and we fall short here of fullfilling this constraint. Rather than looking for a solution of the equation in the class $\mcC^\alpha$ of $\alpha$-H\"older parabolic function, we look for a solution in a restricted class of functions of the form
\begin{equation}  \label{EqPCForm}
u = \widetilde{\sf P}_{u'}Z + u^\sharp,
\end{equation}
for a reference function $Z\in\mcC^\alpha$, to be determined from the noise only and from the equation, with a remainder $u^\sharp\in\mcC^{2\alpha}$ of parabolic H\"older regularity $2\alpha$. Given $Z$, the unknown becomes the pair $(u',u^\sharp)$, with $u'$ in a well-chosen function space. The special paracontrolled form of $u$ allows to make sense of the a priori ill-defined resonant term ${\sf \Pi}(u,\zeta)$, under the assumption that ${\sf \Pi}(Z,\zeta)$ is given as an element of $\mcC^{2\alpha-2}$ -- this Gubinelli, Imkeller and Perkowski's `commutator lemma' \cite{GIP}, Lemma 2.4. We write
$$
\mathscr{L}u = {\sf P}_u\zeta + (2\alpha-2)(u',u^\sharp),
$$
for a function $(2\alpha-2)(u',u^\sharp)$ depending implicitly on $\zeta, Z$ and ${\sf \Pi}(Z,\zeta)$, as a continuous function of all its arguments. From the defining intertwining relation \eqref{EqIntertwining}, the fixed point formulation of equation \eqref{EqPAM} then reads
$$
\widetilde{\sf P}_{u'}Z + u^\sharp = u = \widetilde{\sf P}_u(\mathscr{L}^{-1}\zeta) + \mathscr{L}^{-1}\big((2\alpha-2)(u',u^\sharp)\big),
$$
-- recall we assume for simplicity null initial condition, that is, one has $Z=\mathscr{L}^{-1}(\zeta)$ on the one hand, and 
$$
u'=u= \widetilde{\sf P}_{u'}Z + u^\sharp, \quad u^\sharp = \mathscr{L}^{-1}\big((2\alpha-2)(u',u^\sharp)\big),
$$
on the other hand.

The main feature of the quasilinear setting is the presence of a second order term $Lu$ in the right hand side of the equation. Consider, as a motivation, the model equation
\begin{equation}   \label{EqModelQuasilinear}  \begin{split}
\mathscr{L}u &= u\zeta + uLu   \\
					  &= {\sf P}_u\zeta + {\sf P}_u Lu + \Big({\sf P}_\zeta u + {\sf \Pi}(u,\zeta) + {\sf P}_{Lu}u + {\sf \Pi}(u,Lu)\Big),
\end{split}  \end{equation}
still in the setting where $2/3<\alpha<1$ is close to $1$. As above, one problem is to make sense of the resonant term ${\sf \Pi}(u,Lu)$. This can be done assuming that the term ${\sf \Pi}\big(Z,LZ\big)$ makes sense as an element of the parabolic H\"older space of exponent $2\alpha-2$. With the above well-posedness and regularity assumptions on the resonant term ${\sf \Pi}(Z,\zeta)$, this allows to define the term in parenthesis in the right hand side of \eqref{EqModelQuasilinear} as an element of $\mcC^{2\alpha-2}$. One can see that, for $u$ of paracontrolled form \eqref{EqPCForm}, one has 
$$
{\sf P}_u Lu \simeq {\sf P}_{u'u}LZ,
$$
up to a term in $\mcC^{2\alpha-2}$. A naive fixed point formulation of equation \eqref{EqModelQuasilinear} then reads 
$$
\widetilde{\sf P}_{u'}Z + u^\sharp = \widetilde{\sf P}_u\mathscr{L}^{-1}(\zeta) + \widetilde{\sf P}_{u'u}\mathscr{L}^{-1}(LZ) + (2\alpha)(u',u^\sharp).
$$
Consistency imposes that $Z$ is actually made up of two components $Z=\big(Z^{(1)},Z^{(2)}\big)$, with $Z^{(1)} = \mathscr{L}^{-1}(\zeta)$ and $Z^{(2)} = \mathscr{L}^{-1}(LZ^{(1)})$. The function $u'$ should have as a consequence two components as well, and equation \eqref{EqModelQuasilinear} then rewrites
$$
\sum_{k=1}^2 \widetilde{\sf P}_{u'_k}Z^{(k)} + u^\sharp = \widetilde{\sf P}_u\mathscr{L}^{-1}(\zeta) + \sum_{k=1}^2 \widetilde{\sf P}_{u'_ku}\mathscr{L}^{-1}(LZ^{(k)}) + (2\alpha)(u',u^\sharp),
$$
with terms $\mathscr{L}^{-1}\big({\sf \Pi}(Z^{(i)}, LZ^{(i)})\big)$ inside the remainder $(2\alpha)(\cdots)$ given a priori. The first two terms in the right hand side are taken care of by the $Z^{(1)}$ and $Z^{(2)}$ terms in the left hand side; this is not the case of the term $\widetilde{\sf P}_{u'_2u}\mathscr{L}^{-1}(LZ^{(2)})$ in the right hand side. Consistency then imposes that we actually add a third component to $Z$ and $u'$, to take care of $\widetilde{\sf P}_{u'_2u}\mathscr{L}^{-1}(LZ^{(2)})$. The story then repeats itself, and we are led to consider as a priori form for the solution an infinite paracontrolled expansion
$$
u = \sum_{k\geq 1} \widetilde{\sf P}_{u'_k}Z^{(k)} + u^\sharp,
$$
with $Z^{(k)} = (\mathscr{L}^{-1}L)^{k-1}Z^{(1)}$ for $k>1$, and $Z^{(1)} = \mathscr{L}^{-1}\zeta$. All the $Z^{(k)}$ are elements of $\mcC^\alpha$ here. This infinite dimensional paracontrolled structure is a characteristic feature of the paracontrolled approach of quasilinear singular equations. The convergence of the preceding sum needs to be built in the setting, together with the a priori data of the terms ${\sf \Pi}(Z^{(i)}, LZ^{(i)})$ as elements of $\mcC^{2\alpha-2}$. Anticipating over the results to follow, the reference functions in the paracontrolled expansion of a solution to equation \eqref{EqQgPAM} have the same tree-like structure as the reference functions of a corresponding semilinear equation. This comes from their inductive definition. However, each edge in a `tree' now has a length, corresponding to composing first the operator represented by the edge by the operator $ (\mathscr{L}^{-1}L)^k$, for some $k\geq 0$. This echoes Gerencs\'er and Hairer's work \cite{GerencserHairer}, where each symbol represents an infinite dimensional space. This is the quasilinear effect. The approach works under the quantitative assumption that each a priori term has a natural norm bounded above by a constant multiple of $C^k$, for a constant $C>1$, and $k$ the number of times that the operator $\mathscr{L}^{-1}L$ appears in the formal definition of the term -- the total ``length'' of the tree.

\medskip

We set the scene of paracontrolled calculus in Section \ref{SectionPCScene}, in the form that we need here. Section \ref{SectionProofMainThm} is dedicated to the proof of the well-posedness result in small time for equation \eqref{EqQgPAM}, stated in Theorem \ref{ThmMain}. We give in Appendix \ref{AppendixPCOverview} a bird's eye view on the results from \cite{BB3} on the high order paracontrolled calculus that we use here, while the proofs of a number of new continuity results for operators needed for the study of quasilinear equations are collected in Appendix \ref{AppendixContinuity} and Appendix \ref{SectionExpansionFormula}. 

\bigskip

\noindent \textbf{\textsf{Notations.}} We gather here a number of notations used below.   \vspace{0.15cm}

$\bullet$ \textit{It will be useful sometimes to denote by $(\beta)$ an element of the parabolic H\"older space $\mcC^\beta$ with exponent $\beta$, whose only noticeable feature is its regularity. }

$\bullet$ {\it We denote by $M$ a $3$-dimensional closed Riemannian manifold and set $\mathcal{M} := [0,T]\times M$, for a finite positive time horizon $T$. Given $\alpha\in\bbR$, we denote by $C^\alpha$ the space of $\alpha$-H\"older functions on $M$, defined as the Besov space $B^\alpha_{\infty\infty}$, and write $\mcC^\alpha$ for the parabolic H\"older spaces. We refer the reader to Appendix \ref{AppendixPCOverview} for more information about these spaces.}

\bigskip

\section{Paracontrolled calculus}
\label{SectionPCScene}

One can describe as follows the paracontrolled approach to the study of a generic semilinear singular parabolic PDEs
$$
\mathscr{L} u = f(u,\partial u,\zeta).   
$$
Denote by $P$ the resolution of the free heat equation 
$$
Pu_0:=(\tau,x) \mapsto \big(e^{-\tau L}u_0\big)(x),
$$ 
and recall the intertwining relation \eqref{EqIntertwining} relating $\sf P$ and $\widetilde{\sf P}$.   \vspace{0.2cm}

\begin{enumerate}
   \item[\textsf{\textbf{1{\boldmath $.$}}}] \textsf{\textbf{Paracontrolled ansatz{\boldmath $.$}}} \textit{The irregularity of the noise $\zeta$ dictates the choice of a solution space made up of functions/distributions of the form
   \begin{equation}
   \label{EqAnsatzU}
   u = \sum_{i=1}^{k_0} \widetilde{\sf P}_{u_i}Z_i + u^\sharp,
   \end{equation}
   for reference functions/distributions $Z_i$, of regularity $i\alpha$, that depend formally only on $\zeta$, to be determined later. The order of the expansion is chosen in such a way that $(k_0+1)\alpha+(\alpha-2)>0$. The `derivatives' $u_i$ of $u$ also need to satisfy similar structure equations to a lower order; their derivatives as well, and so on. Denote by $\widehat{u}^\sharp$ the datum of all the remainders in these expansions; they determine entirely this triangular system.}   \vspace{0.15cm}
   
   \item[\textsf{\textbf{2{\boldmath $.$}}}] \textsf{\textbf{Right hand side{\boldmath $.$}}} \textit{Rewrite the right hand side $f(u, \partial u, \zeta)$ of the equation in the canonical form
   \begin{equation}
   \label{EqRHS}
   f\big(u,\partial u, \zeta\big) = \sum_{j=1}^{k_0} {\sf P}_{v_j}Y_j + (\flat)    
   \end{equation}
where $(\flat)$ is a nice remainder and the distributions $Y_j$ depend only on $\zeta$ and the $Z_i$.}   \vspace{0.15cm}
   
   \item[\textsf{\textbf{3{\boldmath $.$}}}]  \textsf{\textbf{Fixed point{\boldmath $.$}}} \textit{The fixed point relation
   \begin{equation*}
   \begin{split}
   u &= Pu_0 + \mathscr{L}^{-1}\big(f(u, \partial u, \zeta)\big)   \\
      &= Pu_0 + \sum_{j=1}^{k_0} \mathscr{L}^{-1}\Big({\sf P}_{v_j}Y_j\Big) + \mathscr{L}^{-1}(\flat)   \\
      &= Pu_0 + \sum_{j=1}^{k_0}\widetilde{\sf P}_{v_j}Z_j + \mathscr{L}^{-1}(\flat),
   \end{split}
   \end{equation*} 
imposes some consistency relations on the choice of the $Z_i = \mathscr{L}^{-1}(Y_i)$ that define them uniquely as functions of $\zeta$, and induces a fixed point relation for $\widehat{u}^\sharp$.}   \vspace{0.15cm}
\end{enumerate}  

Two different questions are addressed in Step \textsf{\textbf{2}}. Making sense of the ill-defined products, characteristic of singular PDEs, and putting the right hand side of the equation in the form \eqref{EqRHS}, for an easy formulation of the fixed point in Step \textbf{\textsf{3}}. One of the main findings of \cite{BB3} is that, at the end of the day, each of these two tasks are dealt with repeating essentially only one operation for each.

Given $\beta\in\bbR$, denote by ${\sf E}^\beta(\cdots)$ a generic (possibly multi-) linear operator that sends continuously $\mcC^\gamma$ into $\mcC^{\beta+\gamma}$, for any $\gamma$ big enough, and such that
\begin{equation}   \label{EqIdentityE}
{\sf E}^\beta(\widetilde{\sf P}_ab,\dots) = a\,{\sf E}^\beta(b,\dots) + {\sf E}^{\beta+\vert b\vert}(a,\dots),
\end{equation}
for all $a\in\mcC^{\vert a\vert}$ and $b\in\mcC^{\vert b\vert}$, with $\vert a\vert, \vert b\vert$ big enough. We say that ${\sf E}$ sends formally $\mcC^\gamma$ into $\mcC^{\beta+\gamma}$ when $\gamma$ is not large enough. A typical example is given by the resonant operator 
$$
{\sf \Pi}(\cdot,c) = {\sf E}^{\vert c\vert}(\cdot),
$$ 
with a fixed argument $c\in\mcC^{\vert c\vert}$; this is part of Gubinelli, Imkeller and Perkowski's important `commutator lemma', Lemma 2.4 in \cite{GIP}. The corrector $\sf C$ from \cite{BB3} and its iterates are $\sf E$-type operators; the operator 
$$
{\sf P}_{\partial a}\partial b = {\sf E}^{-2}(a,b)
$$ 
that appears in the study of the (generalised) (KPZ) equation as well. Another example is
$$
{\sf P}_{Lu}\varepsilon(u,\cdot) = {\sf E}^{-2}\big(u,\varepsilon(u,\cdot)\big).
$$
Applying repeatedly identity \eqref{EqIdentityE} is all we need to investigate the multiplication problem. (The continuity results on the iterated correctors from \cite{BB3} quantify that claim.)

Denote by ${\sf F}^\beta(\cdots)$ a generic (possibly multi-) linear operator that sends continuously $\mcC^\gamma$ into $\mcC^{\beta+\gamma}$, for any $\gamma$, and such that
\begin{equation}   \label{EqIndentityF}
{\sf F}^\beta(\widetilde{\sf P}_ab,\dots) = {\sf P}_a{\sf F}^\beta(b,\dots) + {\sf F}^{\beta+\vert b\vert}(a,\dots),
\end{equation}
for all $a\in\mcC^{\vert a\vert}$ and $b\in\mcC^{\vert b\vert}$, for any $\vert a\vert, \vert b\vert$. Here is a typical example for us
$$
{\sf P}_\zeta u = {\sf F}^{\alpha-2}(u).
$$
Applying repeatedly identity \eqref{EqIndentityF} and continuity results on iterated paraproducts is all we need to put the right hand side in the form \eqref{EqRHS} after all the $\sf E$-operations have been done to analyse the multiplication problems. (The merging operator ${\sf R}^\circ$ from \cite{BB3} is involved here. See Appendix \ref{AppendixPCOverview} for the elements of the high order paracontrolled calculus used in the present work.) With these notations, ${\sf E}^\beta$ and ${\sf F}^\beta$, with no argument, will simply denote elements of $\mcC^\beta$. In those terms, and writing below ${\sf E}^{\beta+|b|}$ for ${\sf E}^{\beta}(b)$, one has for instance
\begin{equation*} \begin{split}
{\sf E}^\beta(\widetilde{\sf P}_ab) &= a\,{\sf E}^{\beta+|b|} + {\sf E}^\beta(a,b)  \\
&= \P_a{\sf E}^{\beta+|b|} + \P_{{\sf E}^{\beta+|b|}}a + {\sf \Pi}\big(a,{\sf E}^{\beta+|b|}\big) + {\sf E}^{\beta+|b|}(a)   \\
&=\P_a{\sf E}^{\beta+|b|} + {\sf F}^{\beta+|b|}(a) + {\sf E}^{\beta+|b|}(a).
\end{split}  \end{equation*}
We can see on this expression that if $a$ itself is given in paracontrolled form $\widetilde{\sf P}_{a_1}a_1$, then we can re-expand the $\sf E$ and $\sf F$ functions of $a$ above. This is the core of the machinery of the high order paracontrolled calculus. We refer the reader to \cite{BB3} for a detailed presentation of the latter. The definitions of the different operators that we use here are recalled in Appendix \ref{AppendixPCOverview}. The quasilinear setting has however two significant features compared to the semilinear setting. Dealing with the second order term $\varepsilon(u, \cdot) Lu$ requires that we work with infinite dimensional paracontrolled system, and one needs to introduce a new corrector together with its iterates to take care of the specific term $\varepsilon(u, \cdot)Lu$.

\bigskip

\subsection{Paracontrolled systems for quasilinear equations}

Fix $0<\alpha<1$. Let an integer $n\geq 1$ be given, together with countable families $\mathscr{T}_1,\dots, \mathscr{T}_n$ of real-valued functions on $[0,T]\times M$, with each $\tau\in \mathscr{T}_i$ of parabolic H\"older regularity $\vert\tau\vert := i\alpha$. Write 
$$
\mathscr{T} := \mathscr{T}_1\cup\cdots\cup \mathscr{T}_n.
$$
A generic finite word with letters in $\mathscr{T}$ will be denoted by $a = (\tau_1,\dots,\tau_k)$, to avoid confusion with the function $\tau_1\cdots\tau_k$, and assigned a \textbf{\textsf{homogeneity}}
$$
\vert a\vert := |\tau_1| + \cdots + |\tau_k|.
$$ 
Define 
$$
\mathscr{A} := \emptyset\cup \Big\{ a = (\tau_1,\dots,\tau_k)\ ;k\geq 1, \ |a| \leq n\alpha\Big\}.
$$
This is the set of words with letters in the alphabet $\mathscr{T}$, and homogeneity no greater than $n\alpha$. This set depends on $n$, which will be fixed in each application. We do no record the dependence of $\mathscr{A}$ on $n$ in the notation. For a word $a = (\tau_1,\dots,\tau_k)$ and $\tau\in \mathscr{T}$, we denote by $a\tau$ the concatenation of $a$ and $\tau$, so $\vert a\tau\vert = \vert a\vert + \vert\tau\vert$. We thus use the symbol $\tau$ both as a function and as a letter in the alphabet $\mathscr{T}$. The setting always makes the meaning of every occurence of a symbol $\tau$ clear; as a rule of thumb $\tau$ is always considered as a letter when it appears in indices. Set $\llparenthesis\emptyset\rrparenthesis := 1$, and for $a=(\tau_1,\ldots,\tau_m)\in\mathscr{A}$, set
$$
\llparenthesis a\rrparenthesis = \|(\tau_1,\ldots,\tau_m)\| := \|\tau_1\|_{\mcC^{|\tau_1|}}\ldots\|\tau_m\|_{\CC^{|\tau_m|}};
$$
this is not a norm. The following definition of a paracontrolled system coincides with the notion used in the study of semilinear singular PDEs, where $\mathscr{T}$ can be chosen to be finite rather than countable.

\medskip

\begin{defn}   \label{pf}
Let $(\beta_a)_{a\in\SA}$ be a family of positive real numbers. A \textbf{\textsf{system paracontrolled by $\mathscr{T}$ at order $n$}} is a family $\widehat{u} = (u_a)_{a\in\SA}$ of parabolic functions such that one has
\begin{equation}   \label{EqPCExpansion}
u_a=\sum_{\tau\in \mathscr{T}; |a\tau|\le n\alpha}\PT_{u_{a\tau}}\tau + u_a^\sharp,
\end{equation}
with $u_a^\sharp \in \mcC^{n\alpha+\beta_a-|a|}$, for all $a\in\mathscr{A}$, and
\begin{equation}   \label{EqConvergenceCondition}
\vvvert\widehat{u}\vvvert := \sum_{b\in\SA}\|u_b^\sharp\|_{\mcC^{n\alpha+\beta_b-|b|}}\llparenthesis b\rrparenthesis < \infty.
\end{equation}
\end{defn}

\medskip

The convergence condition \eqref{EqConvergenceCondition} is always fulfilled in a semilinear setting, where one can work with a finite set $\mathscr{T}$. One proves in Proposition \ref{rc} below that condition \eqref{EqConvergenceCondition} garantees the convergence in a proper space of the (possibly infinite) sum \eqref{EqPCExpansion}. A reasonable choice for the constant $\beta_a$ would be to take them all equal to $\alpha$. This is not a convenient choice from the technical point of view, and all of them will be chosen in the interval $(2/5,\alpha)$ in a particular way explained in Section \ref{SectionProofMainThm} before Theorem \ref{ThmMain}. In particular, they verify $\beta_a>\beta_{a'}$ for any $a,a'\in\SA$ with $a'$ a word containing $a$ as a subword. They play a crucial role in proving that the fixed point formulation of the equation involves a contracting map. We note that all $u_a$ with $\vert a\vert<n\alpha$ are $\mcC^\alpha$, while the $u_a$ with $\vert a\vert=n\alpha$, are elements of $\mcC^{\beta_a}$. Putting together all the contributions from each $\mathscr{T}_i$, each $u_a$ in a paracontrolled system is in particular required to have an expansion of the form
$$
u_a=(\alpha)+(2\alpha)+\ldots+(n\alpha+\beta_a-|a|)
$$
as will be proved in the following propostion. Notice that a paracontrolled system is triangular: The bigger $\vert a\vert$ the lesser we expand $u_a$. Note also that a paracontrolled system is actually determined by the set $\widehat{u}=(u_a^\sharp)_{a\in\mathscr{A}}$ of all remainders in the paracontrolled expansion \eqref{EqPCExpansion}. This motivates that we rewrite the convergence condition \eqref{EqConvergenceCondition} in terms of the remainders only. 

\medskip

\begin{prop} \label{rc}
Let $\widehat{u}$ be a system paracontrolled by $\mathscr{T}$ at order $n$. One has
$$
\sum_{a\in\SA}\|u_a\|_{\CC^{\beta_a}}\llparenthesis a\rrparenthesis \lesssim \vvvert\widehat{u}\vvvert.
$$
This implies in particular 
$$
\|u_a\|_{\CC^{\beta_a}}\lesssim\vvvert\widehat{u}\vvvert, \quad \forall\,a\in\mathscr{A}.
$$
\end{prop}

\medskip

\begin{Dem}
Given any $a\in\SA$, we have by a finite induction
\begin{align*}
\|u_a\|_{\CC^{\beta_a}}&\lesssim\sum_{\tau\in\mathscr{T};|a\tau|\le n\alpha}\|u_{a\tau}\|_{\CC^{\beta_{a\tau}}}\|\tau\|+\|u_a^\sharp\|_{\CC^{\beta_a}}\\
&\lesssim\sum_{b\in\SA;|ab|\le n\alpha}\|u_{ab}^\sharp\|_{\CC^{\beta_{ab}}}\llparenthesis b\rrparenthesis.
\end{align*}
This yields
$$
\sum_{a\in\SA}\|u_a\|_{\CC^{\beta_a}}\llparenthesis a\rrparenthesis \lesssim \sum_{a\in\SA}\|u_a^\sharp\|_{\CC^{\beta_a}}\llparenthesis a\rrparenthesis \lesssim \vvvert\widehat{u}\vvvert.
$$
\end{Dem}

\bigskip

\subsection{Additional correctors}
\label{SectionAdditionalCorrectors}

The formulation of the quasilinear equation \eqref{EqQgPAM} in the semilinear-like form \eqref{EqReformulation} involves the second order term $\varepsilon(u,\cdot)Lu$, specific to the quasilinear setting. Writing
\begin{equation}   \label{EqDecompositionEpsilon}
\varepsilon(u,\cdot) Lu = {\sf P}_{\varepsilon(u,\cdot)}Lu + {\sf P}_{Lu}\varepsilon(u,\cdot) + {\sf \Pi}\big(\varepsilon(u,\cdot),Lu\big),
\end{equation}
the $\sf E$-type operators 
$$
{\sf P}_{La}b, \quad {\sf \Pi}(La,b)
$$
that appear in the last two terms of the right hand side of identity \eqref{EqDecompositionEpsilon} happen to be of the same type as the resonant operator 
$(a,b)\mapsto {\sf \Pi}(a,b)$. Their analysis is thus similar to what was done in \cite{BB3} for the resonant operator via the introduction of the corrector $\sf C$ and its iterates. The $\sf F$-type operator 
$$
{\sf P}_a Lb
$$
that appears in the first term of the right hand side of \eqref{EqDecompositionEpsilon} does not show up in the study of semilinear singular PDEs and requires a specific treatment. We state here a number of continuity results whose proofs are given in Appendix \ref{AppendixContinuity}; all the proofs are variations on the pattern of proofs of continuity results from \cite{BB3}. Given that the technical setting of \cite{BB2, BB3} is likely not to be familiar to most readers, we also give in this section the proofs of some of the statements in the time-independent model setting of the flat torus. The paraproduct and resonant operators 
$$
P^0_ab := \sum_{i<j-1}\Delta_i(a)\Delta_j(b), \quad \Pi^0(a,b) := \sum_{\vert i-j\vert\leq 1}\Delta_i(a)\Delta_j(b),
$$
are then defined classically in terms of Fourier projectors $\Delta_k$. We refer the reader to \cite{BCD} for the basics on Littlewood-Paley decomposition and paraproduct and resonant operators in that setting.

The continuity results from this section are all we need in addition to the results of \cite{BB3} to study equation \eqref{EqQgPAM}, and more generally a whole class of quasilinear singular PDEs.

\bigskip

\noindent \textbf{\textsf{2.2.1 Operator ${\sf P}_a Lb$.}} We define the operator 
$$
{\sf L}(a,b) := L\PT_ab - \P_aLb.
$$
Continuity results on this operator allow to get an expansion for $Lu$ of the form
$$
Lu = \sum {\sf P}_{u'_\tau}(L\tau) + (4\alpha-2),
$$
for some $u'_\tau$, from a paracontrolled expansion for $u$. A paracontrolled expansion for a term of the form ${\sf P}_a(Lu)$ can then be obtained. We also define the refined operator
$$
{\sf L}_{(1)}(a,b) := L\PT_ab-\P_aLb - \sum_{i=1}^\ell \P_{d(\overline{u}_0)^{-1}V_ia}^{(i)}Lb
$$
to deal with arguments $a$ in ${\sf L}(a,b)$ with regularity exponent greater than $1$. The operators $\P^{(i)}$ are defined by for any $e$ in the parabolic space $\CM$ by
$$
\left(\P_a^{(i)}b\right)(e):=\int_{e',e''\in\CM}K(e;e',e'')a(e')\left(\PT_{\delta_i(\cdot,e')}b\right)(e'')\,\nu(\drm e')\nu(\drm e'')
$$
with $K$ the kernel of the bilinear operator $(a,b)\mapsto\P_ab$. See Appendix \ref{AppendixPCOverview} for notations and details on the parabolic setting. The following theorem is proved in Theorem \ref{ThmLV} in Appendix \ref{AppendixContinuity}.

\medskip

\begin{thm}   \label{ThmContinuitySfL}
\begin{itemize}
	\item Let $\alpha\in(0,1)$ and $\beta\in(-3,3)$, be such that $\alpha+\beta<3$, and $\alpha+\beta-2\in(-3,3)$. Then the operator ${\sf L}$ has a natural extension as a continuous operator from $\mcC^{\alpha}\times\mcC^{\beta}$ into $\mcC^{\alpha+\beta-2}$.   \vspace{0.1cm}
	
	\item Let $\alpha_1,\alpha_2\in(0,1)$ and $\beta\in(-3,3)$ such that $\alpha_1+\beta<3$ and $\alpha_1+\alpha_2+\beta-2\in(-3,3)$. Then the iterated operator 
	$$
	{\sf L}\big((a_1,a_2),b\big) := {\sf L}\big(\P_{a_1}a_2,b\big) - \P_{a_1}{\sf L}(a_2,b)
	$$ 
	has a natural extension as a continuous operator from $\mcC^{\alpha_1}\times\mcC^{\alpha_2}\times\mcC^\beta$ into $\mcC^{\alpha_1+\alpha_2+\beta-2}$.   \vspace{0.1cm}
	
	\item Let $\alpha_1,\alpha_2,\alpha_3\in(0,1)$ and $\beta\in(-3,3)$ such that $\alpha_1+\alpha_2+\beta<3$, $\alpha_2+\beta<3$ and $\alpha_1+\alpha_2+\beta-2\in(-3,3)$. Then the iterated operator 
	$$
	{\sf L}\Big(\big((a_1,a_2),a_3\big), b\big) := {\sf L}\big((\P_{a_1}a_2,a_3), b\big) - \P_{a_1}{\sf L}\big((a_2,a_3), b\big)
	$$ 
	has a natural extension as a continuous operator from $\mcC^{\alpha_1}\times\mcC^{\alpha_2}\times\mcC^{\alpha_3}\times\mcC^\beta$ into $\mcC^{\alpha_1+\alpha_2+\alpha_3+\beta-2}$.    \vspace{0.1cm}

    \item Let $\alpha\in(1,2)$ and $\beta\in(-3,3)$, be such that $\alpha+\beta<3$, and $(\alpha+\beta-2)\in(-3,3)$. Then the operator ${\DL_{(1)}}$ has a natural extension as a continuous operator from $\CC^{\alpha}\times\CC^{\beta}$ into $\CC^{\alpha+\beta-2}$.
\end{itemize}
\end{thm}

\bigskip

\noindent \textbf{\textsf{2.2.2 Operators ${\sf P}_{La} b$ and ${\sf \Pi}(La,b)$.}} These two operators of $\sf E$-type are defined by similar formulas as the resonant operator, in terms of the parabolic approximation operators $\mcQ_t$ from \cite{BB3}. It is thus natural that they satisfy expansion rules similar to the expansion rules satisfied by the resonant operator. Introduce for that purpose the operators
\begin{align*}
{\sf C}_L^<\Big((a_1,a_2), b\Big) &:= \P_{L\PT_{a_1}a_2}b - a_1\P_{La_2}b,   \\
{\sf C}_L^>\Big(a, (b_1,b_2)\Big) &:= \P_{La}\Big(\PT_{b_1}b_2\Big) - b_1\P_{La}b_2,   \\
{\sf C}_L\Big((a_1,a_2), b\Big) &:= \PI\Big(L\PT_{a_1}a_2 ,b\Big) - a_1\PI\Big(La_2, b\Big).
\end{align*}

We choose the notation $<$ in the exponent of ${\sf C}_L^<$ to emphasize that the paraproduct term is in the low `frequency' part of the operator, while it is in the high `frequency' part in ${\sf C}_L^>$. The following theorem is proved here in the time-independent model setting of the flat torus; see Theorem \ref{ThmCL} in Appendix \ref{AppendixContinuity} for the proof.

\medskip

\begin{thm}   \label{ThmContinuityCL}
\begin{itemize}
	\item Let $\alpha_1\in(0,1)$ and $\alpha_2,\beta\in(-3,3)$ such that $\alpha_1+\alpha_2\in(-3,3)$. If
	\begin{equation}   \label{EqConditionsExponents}
	\alpha_2+\beta-2<0\quad\text{and}\quad\alpha_1+\alpha_2+\beta-2>0
	\end{equation}
	then the operators ${\sf C}_L^<$ and ${\sf C}_L$ have natural extensions as continuous operators from $\mcC^{\alpha_1}\times\mcC^{\alpha_2}\times\mcC^\beta$ into $\mcC^{\alpha_1+\alpha_2+\beta-2}$.   \vspace{0.1cm}
	
	\item Let $\beta_1\in(0,1)$ and $\alpha,\beta_2\in(-3,3)$ such that $\beta_1+\beta_2\in(-3,3)$. If
	$$
	\alpha+\beta_2-2<0\quad\text{and}\quad\alpha+\beta_1+\beta_2-2>0
	$$
	then the operator ${\sf C}_L^>$ has a natural extension as a continuous operator from $\mcC^\alpha\times\mcC^{\beta_1}\times\mcC^{\beta_2}$ into $\mcC^{\alpha+\beta_1+\beta_2-2}$.
\end{itemize}
\end{thm}

\medskip

\begin{Dem}
Write $\Delta$ for the usual Laplacian on the flat torus. 

\ssk

$\bullet$ Set
$$
C_\Delta^0(a_1,a_2,b) := \Pi^0(\Delta P^0_{a_1}a_2, v) - a_1\Pi^0(\Delta a_2, b).
$$
We prove that for $\alpha_1,\alpha_2$ and $\beta$ such that inequalities \eqref{EqConditionsExponents} hold true, the operator $C_\Delta^0$ is continuous from $C^{\alpha_1}\times C^{\alpha_2}\times C^\beta$ into $C^{\alpha_1+\alpha_2+\beta-2}$. We have
$$
C_\Delta^0(a_1,a_2,b) = \sum_{|i-j|<1} \Delta_i\left(P^0_{a_1}a_2\right)\Delta_j(b) - a_1\Delta_i(a_2)\Delta_j(b).
$$
Setting 
$$
\varepsilon_i := \Delta_i\big(\Delta P^0_{a_1}a_2\big) - a_1\Delta_i(\Delta a_2),
$$ 
we have
$$
C_\Delta^0(a,b,c) = \sum_{|i-j|<1} \varepsilon_i\,\Delta_j(b).
$$
As in the proof of the estimate for the classic corrector $\DC$, one sees that one has
$$
\|\Delta_k\varepsilon_i\|_{L^\infty}\lesssim 2^{2i}\,2^{-i\alpha_2}2^{-\max(i,k)\alpha_1}\,\|a_1\|_{C^{\alpha_1}}\|a_2\|_{C^{\alpha_2}};
$$
the factor $2^{2i}$ comes from the $\Delta$ operator. Writing
\begin{align*}
\Delta_k\left(C_\Delta^0(a_1,a_2,b)\right) &= \sum_{|i-j|\le1} \Delta_k\big(\varepsilon_i\,\Delta_j(b)\big)   \\
&= \sum_{\underset{|i-j|\le1}{j<k-2}} \Delta_k(\varepsilon_i)\Delta_j(b) + \sum_{\underset{|i-j|\le1}{k<j-2}} \Delta_k\big(\Delta_i(\varepsilon_i)\Delta_j(b)\big) + \sum_{\underset{|i-j|\le1}{|k-j|\le1}} \Delta_k\big(S_i(\varepsilon_i)\Delta_j(b)\big),
\end{align*}
we see that
\begin{small}\begin{align*}
\Big\|\Delta_k\left(C_\Delta^0(a_1,a_2,b)\right)\Big\|_{L^\infty} &\lesssim \left\{\sum_{i<k-2}2^{-i(\alpha_2+\beta-2)}2^{-k\alpha_1}+\sum_{k<i-2}2^{-i(\alpha_1+\alpha_2+\beta-2)}\right.   \\
&\quad\quad\left.+\sum_{|i-k|\le1}2^{-i(\alpha_1+\alpha_2+\beta-2)}\right\}\,\|a_1\|_{C^{\alpha_1}}\|a_2\|_{C^{\alpha_2}}\|b\|_{C^\beta}   \\
&\lesssim 2^{-k(\alpha_1+\alpha_2+\beta-2)}\,\|a_1\|_{C^{\alpha_1}}\|a_2\|_{C^{\alpha_2}}\|b\|_{C^\beta}
\end{align*}\end{small}
using that $(\alpha_2+\beta-2)<0$ and $(\alpha_1+\alpha_2+\beta-2)>0$.

\medskip

$\bullet$ Set now
\begin{equation*} \begin{split}
C_\Delta^{<,0}(a_1,a_2,b) :=& P^0_{\Delta P^0_{a_1}a_2}b - a_1P^0_{\Delta a_2}b   \\
=& \sum_{i<j-2} \varepsilon_i\,\Delta_j(b).
\end{split} \end{equation*}
We prove that for $\alpha_1,\alpha_2$ and $\beta$ such that inequalities \eqref{EqConditionsExponents} hold true, $C_\Delta^{<,0}$ is continuous from $C^{\alpha_1}\times C^{\alpha_2}\times C^\beta$ into $C^{\alpha_1+\alpha_2+\beta-2}$.
This can be seen by writing
\begin{align*}
\Delta_k\left(\DC_L^{<,0}(a_1,a_2,b)\right) = \sum_{\underset{i<j-2}{j<k-2}} \Delta_k(\varepsilon_i)\Delta_j(b) + \sum_{\underset{i<j-2}{k<j-2}}\Delta_k\big(\Delta_i(\varepsilon_i)\Delta_j(b)\big) + \sum_{\underset{i<j-2}{|k-j|\le1}}\Delta_k\big(S_i(\varepsilon_i)\Delta_j(b)\big),
\end{align*}
from which one sees that
\begin{small}\begin{align*}
\left\|\Delta_k\left\{\DC_L^{<,0}(a_1,a_2,b)\right)\right\|_{L^\infty} &\lesssim \left\{\sum_{j<k-2}2^{-i(\alpha_2+\beta-2)}2^{-k\alpha_1}+\sum_{k<j-2}2^{-i(\alpha_1+\alpha_2+\beta-2)}\right.   \\
&\quad\quad\left.+\sum_{|j-k|\le1}2^{-i(\alpha_1+\alpha_2+\beta-2)}\right\}\,\|a_1\|_{C^{\alpha_1}}\|a_2\|_{C^{\alpha_2}}\|b\|_{C^\beta}   \\
&\lesssim 2^{-k(\alpha_1+\alpha_2+\beta-2)}\,\|a_1\|_{C^{\alpha_1}}\|a_2\|_{C^{\alpha_2}}\|b\|_{C^\beta}.
\end{align*}\end{small}
\end{Dem}

\medskip

We use this continuity result under the form of the $\sf E$-type identity
\begin{equation} \label{EqIdentity1}
{\sf P}_{L\widetilde{\sf P}_{a_1}a_2}b = {\sf E}^{-2}({\sf P}_{a_1}a_2,b) = a_1{\sf E}^{-2}(a_2,b) + {\sf E}^{-2+\vert a_2\vert}(a_1,b),
\end{equation}
or its analogue with $b$ given by a paraproduct; Theorem \ref{ThmContinuityCL} justifies fully this identity in the regime $\alpha_1\in(0,1)$ or $\beta_1\in(0,1)$. In a setting where $a_2$ and $b$ play the role of data, one rewrites identity \eqref{EqIdentity1} as
$$
{\sf P}_{L\widetilde{\sf P}_{a_1}a_2}b = {\sf E}^{-2}({\sf P}_{a_1}a_2,b) = {\sf P}_{a_1}{\sf E}^{-2+\vert a_2\vert+\vert b\vert} + {\sf F}^{-2+\vert a_2\vert+\vert b\vert}(a_1) + {\sf E}^{-2+\vert a_2\vert+\vert b\vert}(a_1).
$$
We also have continuity estimates on iterated correctors, as in \cite{BB3}. Given the proof of Theorem \ref{ThmContinuityCL} given in Appendix \ref{AppendixContinuity}, it will be clear to the reader that their statements and proofs are identical to what is done in \cite{BB3} for the iterated correctors, see Section 3.1.3 therein. We leave their statements and proofs to the reader. The continuity results from Theorem \ref{ThmContinuityCL} take profit only from the H\"older regularity of $a_1$ or $b_1$, for any regularity exponent in $(0,1]$. As in the semilinear case, we need to introduce refined correctors to refine the estimates if $a_1$ or $b_1$ are $\alpha_1$ or $\beta_1$-Lipscthiz, with $\alpha_1$ or $\beta_1$ strictly greater than $1$. We set for that purpose, for a generic spacetime point $e$,
\begin{align*}
{\sf C}_{L,(1)}^<\Big(a_1,a_2,b\Big)(e) &:= {\sf C}_L^<(a_1,a_2,b)(e) - d\big(\overline{u_0}(e)\big)^{-1} \sum_{i=1}^\ell  (V_ia_1)(e) \, \left(\P_{L\PT_{\delta_i(e,\cdot)}a_2}b\right)(e),   \\
{\sf C}_{L,(1)}^>\Big(a,b_1,b_2\Big)(e) &:= {\sf C}_L^>(a,b_1,b_2)(e) - d\big(\overline{u_0}(e)\big)^{-1} \sum_{i=1}^\ell (V_ib_1)(e) \, \left(\P_{La}\PT_{\delta_i(e,\cdot)}b_2\right)(e), 
\end{align*}
\begin{align*}
{\sf C}_{L,(1)}\Big(a_1,a_2,b\Big)(e) &:= {\sf C}_L(a_1,a_2,b)(e) - d\big(\overline{u_0}(e)\big)^{-1}  \sum_{i=1}^\ell (V_ia_1)(e) \, \PI\Big(L\PT_{\delta_i(e,\cdot)}a_2,b\Big)(e),
\end{align*}
where the functions $\delta_i$ are defined in Appendix \ref{AppendixContinuity}. Keep in mind right now that in the setting of the flat torus $d\big(\overline{u_0}(e)\big)^{-1}V_i = \partial_i$, the partial derivative in the $i^\textrm{th}$ space direction, and $\delta_i(e,e') = d\big(\overline{u_0}(x)\big)^{1/2}(x_i - x'_i)$, for spacetime points $e=(t,x)$ and $e'=(t',x')$. The following theorem is also proved here in the time-independent model setting of the flat torus; see Theorem \ref{ThmCL1} in Appendix \ref{AppendixContinuity} for the proof.

\medskip

\begin{thm}   \label{ThmContinuityRefinedCL}
\begin{itemize}
	\item Let $\alpha_1\in(1,2)$ and $\alpha_2,\beta\in(-3,3)$ such that $\alpha_1+\alpha_2\in(-3,3)$. If
	$$
	\alpha_2+\beta-2<0\quad\text{and}\quad\alpha_1+\alpha_2+\beta-2>0
	$$
	then the operators ${\sf C}_{L,(1)}^<$ and ${\sf C}_{L,(1)}$ have natural extensions as continuous operators from $\mcC^{\alpha_1}\times\mcC^{\alpha_2}\times\mcC^\beta$ into $\mcC^{\alpha_1+\alpha_2+\beta-2}$.   \vspace{0.1cm}
	
	\item Let $\beta_1\in(1,2)$ and $\alpha,\beta_2\in(-3,3)$ such that $\beta_1+\beta_2\in(-3,3)$. If
	$$
	\alpha+\beta_2-2<0\quad\text{and}\quad\alpha+\beta_1+\beta_2-2>0
	$$
	then the operator ${\sf C}_{L,(1)}^>$ has a natural extension as a continuous operator from $\mcC^\alpha\times\mcC^{\beta_1}\times\mcC^{\beta_2}$ into $\mcC^{\alpha+\beta_1+\beta_2-2}$.
\end{itemize}
\end{thm}

\medskip

In terms of $\sf E$-type identities, this continuity result rewrites under the form
$$
{\sf P}_{L\widetilde{\sf P}_{a_1}a_2}b = {\sf E}^{-2}({\sf P}_{a_1}a_2,b) = a_1{\sf E}^{-2}(a_2,b) + \sum_{i=1}^\ell d_0^{-1} (V_ia_1) {\sf E}_i^{-1}(a_2,b) +  {\sf E}^{-2+\vert a_2\vert}(a_1,b),
$$
with $d_0 := d(\overline{u_0})$, and similar expressions for ${\sf C}^<_{L,(1)}$ and ${\sf C}^>_{L,(1)}$. This identity holds here in the regime $1<\vert a_1\vert<2$. In a setting where $a_2$ and $b$ play the role of data, one rewrites the preceding identity as
$$
{\sf P}_{L\widetilde{\sf P}_{a_1}a_2}b = {\sf P}_{a_1}{\sf E}^{-2+\vert a_2\vert+\vert b\vert} + {\sf F}^{-2+\vert a_2\vert+\vert b\vert}(a_1) + {\sf E}^{-2+\vert a_2\vert+\vert b\vert}(a_1).
$$
This identity takes here the same form as in the regime $0<a_1<1$. This is the form that we use in the computations.

\medskip

Theorem \ref{ThmContinuitySfL}, Theorem \ref{ThmContinuityCL} and Theorem \ref{ThmContinuityRefinedCL} take care of the specific features of quasilinear equations, compared to their semilinear analogue. Formulation \eqref{EqReformulation} also involve a term $a_i(u,\cdot)V_iu$ that can appear in a semilinear setting as well, and the function $\varepsilon(u,\cdot)$. The last two paragraphs of this section state the results that we need about them.

\bigskip

\noindent \textbf{\textsf{2.2.3 Dealing with the term $a_i(u,\cdot)V_iu$.}} We have the following continuity results for the operators
\begin{align*}
{\sf C}_{V_i}^<\Big(a_1,a_2,b\Big) &:= \P_{V_i\PT_{a_1}a_2}b - a_1\P_{V_ia_2}b,   \\
{\sf C}_{V_i}^>\Big(a,b_1,b_1\Big) &:= \P_{V_ia}\Big(\PT_{b_1}b_2\Big) - b_1\P_{V_ia}b_2,   \\
{\sf C}_{V_i}\Big(a_1,a_2,b\Big) &:= \PI\Big(V_i\PT_{a_1}a_2,b\Big) - a_1\PI\Big(V_ia_2,b\Big);
\end{align*}
see Theorem \ref{ThmCL} in Appendix \ref{AppendixContinuity} for the proof. 

\medskip

\begin{thm}   \label{ThmVi}
\begin{itemize}
	\item Let $\alpha_1\in(0,1)$ and $\alpha_2,\beta\in(-3,3)$ such that $\alpha_1+\alpha_2\in(-3,3)$. If
	\begin{equation} \label{EqIdentity2}
	\alpha_2+\beta-1<0\quad\text{and}\quad\alpha_1+\alpha_2+\beta-1>0
	\end{equation}
	then the operators ${\sf C}_{V_i}^<$ and ${\sf C}_{V_i}$ have natural extensions as continuous operators from $\mcC^{\alpha_1}\times\mcC^{\alpha_2}\times\mcC^\beta$ into $\mcC^{\alpha_1+\alpha_2+\beta-1}$.   \vspace{0.1cm}
	
	\item Let $\beta_1\in(0,1)$ and $\alpha,\beta_2\in(-3,3)$ such that $\beta_1+\beta_2\in(-3,3)$. If
	$$
	\alpha+\beta_2-1<0\quad\text{and}\quad\alpha+\beta_1+\beta_2-1>0
	$$
	then the operator ${\sf C}_{V_i}^>$ has a natural extension as a continuous operator from $\mcC^\alpha\times\mcC^{\beta_1}\times\mcC^{\beta_2}$ into $\mcC^{\alpha+\beta_1+\beta_2-1}$.
\end{itemize}
\end{thm}

\medskip

\begin{Dem}
We prove here this continuity result for a simplified version of the operator $\DC_V$ in the time-independent case of the flat torus, with the constant vector field $\partial_1$ in the role of $V_i$; we refer the reader to Appendix \ref{AppendixContinuity} for the proof of Theorem \ref{ThmVi} in the general setting. Set
$$
C_{\partial_1}^0(a,b,c) := \Pi^0\big(\partial_1 P^0_ab, c\big) - a \Pi^0(\partial_1 b, c).
$$
We prove that for $\alpha,\beta$ and $\gamma$ such that inequalities \eqref{EqIdentity2} hold true, the operator $C_{\partial_1}^0$ is continuous from $C^\alpha\times C^\beta\times C^\gamma$ into $C^{\alpha+\beta+\gamma-2}$. Using that $\Delta_i(\partial_1 f)\simeq O(2^i)\Delta_i(f)$, for a function $O(2^i)$ with uniform norm of order $2^i$, we have
$$
C_{\partial_1}^0(a,b,c)\simeq\sum_{|i-j|<1}O(2^i)\Delta_i\left(\PI^0_ab\right)\Delta_j(c) - aO(2^i)\Delta_i(b)\Delta_j(c),
$$
so
$$
C_{\partial_1}^0(a,b,c) = \sum_{|i-j|<1} O(2^i)\varepsilon_i\Delta_j(c).
$$
The same computations as above then yield the estimate
$$
\left\|\Delta_k\left(\DC_V^0(a,b,c)\right)\right\|_{L^\infty}\lesssim 2^{-k(\alpha+\beta+\gamma-1)}\|a\|_{C^\alpha}\|b\|_{C^\beta}\|c\|_{C^\gamma}.
$$
\end{Dem}

\medskip

Theorem \ref{ThmVi} justifies that we summarize the above continuity statement under the following $\sf E$-type identity
$$
\P_{V_i\PT_{a_1}a_2}b = a_1\P_{V_ia_2}b + {\sf E}^{\vert a_2\vert-1}(a_1,b),
$$
with similar identities satisfied by the expressions $\P_{V_ia}\Big(\PT_{b_1}b_2\Big)$ and $\PI\Big(V_i\PT_{a_1}a_2,b\Big)$. In a setting where $a_2$ and $b$ play the role of data, one rewrites the preceding identity as
$$
\P_{V_i\PT_{a_1}a_2}b = {\sf P}_{a_1}{\sf E}^{\vert a_2\vert+\vert b\vert-1} + {\sf E}^{\vert a_2\vert+\vert b\vert-1}(a_1).
$$
This is the form under which we use Theorem \ref{ThmVi} in computations.

\medskip

Associate with each vector field $V_i$ the operator
$$
\DV_i (a,b) := V_i\big(\widetilde{\sf P}_ab\big) - {\sf P}_a(V_ib).
$$
We prove the theorem here in the time-independent model setting of the flat torus; see Theorem \ref{ThmLV} in Appendix \ref{AppendixContinuity} for the proof.

\medskip

\begin{thm}   \label{ThmVi2}
\begin{itemize}
	\item Let $\alpha,\beta\in(-3,3)$ such that $\alpha+\beta-1\in(-3,3)$. Then the operator $\DV_i$ has a natural extension as a continuous operator from $\mcC^\alpha\times\mcC^\beta$ to $\mcC^{\alpha+\beta-1}$.   \vspace{0.1cm}
	
	\item Let $\alpha_1,\alpha_2\in(0,1)$ and $\beta\in(-3,3)$ such that $\alpha_1+\beta<3$ and $\alpha_1+\alpha_2+\beta-1\in(-3,3)$. Then the iterated operator 
	$$
	\DV_i((a_1,a_2),b):=\DV_i(\widetilde{\P}_{a_1}a_2,b) - \P_{a_1}\DV_i(a_2,b)
	$$ 
	has a natural extension as a continuous operator from $\mcC^{\alpha_1}\times\mcC^{\alpha_2}\times\mcC^\beta$ to $\mcC^{\alpha_1+\alpha_2+\beta-1}$.
\end{itemize}
\end{thm}

\bigskip

\noindent \textbf{\textsf{2.2.4 Paracontrolled expansion of $\varepsilon(u,\cdot)$.}} Finally, we have the following variation on the high order paracontrolled expansion formula from \cite{BB3}, Theorem 4 therein. 

\medskip

\begin{thm}  \label{ThmParacontrolledExpansion}
Let $f:\IR\to\IR$, be a $C_b^4$ function and let $u$ and $v$ be respectively $\mcC^\alpha$ and $\mcC^{4\alpha}$ functions on $[0,T]\times M$, with $\alpha\in(0,1)$. Then we have
\begin{align*}
f(u)v&=\P_{f'(u)v}u + \frac{1}{2}\Big\{\P_{f^{(2)}(u)v}u^2-2\P_{f^{(2)}(u)uv}u\Big\}   \\
&\quad+ \frac{1}{3!}\Big\{\P_{f^{(3)}(u)v}u^3-3\P_{f^{(3)}(u)uv}u^2+3\P_{f^{(3)}(u)u^2v}u\Big\} + (\sharp),
\end{align*}
for a remainder $(\sharp)\in\mcC^{4\alpha}$.
\end{thm}

\medskip

The proof of this statement is given in Appendix \ref{SectionExpansionFormula}.

\bigskip

\section{Quasilinear generalised (PAM) equation}
\label{SectionProofMainThm}

We use the generic three step process from Section \ref{SectionPCScene} to solve the quasilinear generalised (PAM) equation \eqref{EqReformulation}. 

\medskip

\textbf{\textsf{Step 1.}} We have $2/5<\alpha<1/2$, so we choose to work with a third order paracontrolled expansion, to have a remainder term $u^\sharp$ in the paracontrolled expansion of $u$ for which the product of $u^\sharp\in\mcC^{4\alpha}$ with any distribution of H\"older regularity $\alpha-2$ is well-defined.

\medskip

\textbf{\textsf{Step 2.}} We use continuity results for correctors, commutators and their iterates, to put the right hand side of equation \eqref{EqReformulation} in the canonical form \eqref{EqRHS}. Recall $d_0=d(\overline{u_0})$. Indices $a,b,c$ below are in $\mathscr{A}$, while $\tau\in\mathscr{T}$.

\medskip

\begin{prop}   \label{PropDecomposition}
Assume we are given a system $(u_a)_{a\in\mathscr{A}}$ paracontrolled by a family $\mathscr{T}$ at order $3$. Then
\begin{equation}   \label{EqDecomposition}  \begin{split}
f(u)\zeta &+ \varepsilon(u, \cdot)Lu + \sum_{i=1}^\ell a_i(u,\cdot)V_iu   \\
&=\P_{f(u)}\zeta + \sum_{|a|\le2\alpha} \P_{f'(u)u_a}\zeta_a^{(1)} + \sum_{|ab|\le2\alpha} \P_{f^{(2)}(u)u_au_b}\zeta_{ab}^{(1)}   \\ 
&\quad+ \sum_{\tau\in\mathscr{T}} {\sf P}_{\varepsilon(u,\cdot)u_\tau}L\tau + \sum_{|a|\le 3\alpha; a\in\mathscr{A}\backslash\mathscr{T}} \P_{\varepsilon(u)u_a}\zeta_a^{(2)}   \\
&\quad+ \sum_{|ab|\le3\alpha}\P_{d_0^{-1}d'(u)u_au_b}\zeta_{ab}^{(2)} + \sum_{|abc|\le3\alpha} \P_{d_0^{-1}d^{(2)}(u)u_au_bu_c} \zeta_{abc}^{(2)}   \\
&\quad+ \sum_{|\tau|=\alpha; 1\leq j\leq \ell} \P_{a_j(u,\cdot)u_\tau} \zeta_{j,\tau} + (\sharp),
\end{split} \end{equation}
for distributions $\zeta^{(1)}_e, \zeta^{(2)}_e, \zeta_{j,\tau}$ that depend only on $\zeta$ and $\mathscr{T}$, with $\zeta^{(1)}_e$ of regularity $\vert e\vert+\alpha-2$, with $\zeta^{(2)}_e$ of regularity $\vert e\vert-2$ and $\zeta_{j,\tau}$ of regularity $\vert\tau\vert-1$, for $e\in\mathscr{A}, \tau\in\mathscr{T}$ and $1\leq j\leq \ell$. The remainder $(\sharp)$ is an element of $\mcC^{4\alpha-2}$.
\end{prop}

\medskip

As always in the analytic part of the study of a singular PDE, one needs to assume that the distributions $\zeta_e^{(1)}, \zeta_e^{(2)}, \zeta_{j,\tau}$ are given off-line. The remainder term $(4\alpha-2)$ also involves off-line data. The point with \textit{stochastic} singular PDEs is that one can construct the data by probabilistic means; this is what renormalisation is about. It comes as a by-product of the proof that the remainder is the sum of a term of regularity $4\alpha-2$ involving off-line data and a term of regularity $5\alpha-2$ that is a continuous function of the paracontrolled system $(u_a)_{a\in\mathscr{A}}$ and all the off-line data.

\medskip

\begin{Dem}
Below, we check the convergence of all implicit infinite sums  of $\mathscr{T}$ using the convergence condition \eqref{EqConvergenceCondition} in the definition of a paracontrolled system; we do not do that explicitly each time. Recall we denote by $(\beta)$ an element of the parabolic H\"older space $\mcC^\beta$ with regularity exponent any $\beta\in\bbR$, whose only noticeable feature is its regularity. Its expression may change from line to line. Recall also from Appendix \ref{AppendixPCOverview} the definition of the operator
$$
{\sf R}^\circ(a,b,c) = {\sf P}_a{\sf P}_bc - {\sf P}_{ab}c,
$$
its continuity and expansion properties. To shorten notations, we sometimes use implicit summation on repeated indices.   \vspace{0.15cm}

\begin{itemize}
	\item The term $f(u)\zeta$ is the same as in the semilinear (gPAM) equation so its decomposition is given by proposition $17$ of \cite{BB3}, that is
	$$
	f(u)\zeta = \P_{f(u)}\zeta + \sum_{|a|\le2\alpha}\P_{f'(u)u_a} \zeta_a^{(1)} + \sum_{|ab|\le2\alpha}\P_{f^{(2)}(u)u_au_b} \zeta_{ab}^{(1)}+ (4\alpha-2).
	$$

	\item For the term $\P_{\varepsilon(u,\cdot)}Lu$, first, we have from Theorem \ref{ThmContinuitySfL}
	\begin{align*}
	Lu &= L\PT_{u_\tau}\tau+(4\alpha-2)   \\
	    &= \P_{u_\tau}L\tau + \DL(u_\tau,\tau) + (4\alpha-2)   \\
		&= \P_{u_\tau}L\tau + \P_{u_{\tau\Mj}}\DL(\Mj,\tau) + \DL((u_{\tau\Mj},\Mj),\tau)+(4\alpha-2)   \\
		&= \P_{u_\tau}L\tau + \P_{u_{\tau\Mj}}\DL(\Mj,\tau) + \P_{u_{\tau\Mj\Mk}}\DL((\Mk,\Mj),\tau) + (4\alpha-2).
	\end{align*}
	One takes care of remainder terms in the expansions of the $u_\tau$'s with $\vert\tau\vert=\alpha$, in the expression ${\sf L}(u_\tau,\tau)$, using the operator ${\sf L}_{(1)}$. Write the above expression under the form
	$$
	Lu =: {\sf P}_{u_a}\xi^{(2)}_a + (4\alpha-2),
	$$
	with $\xi^{(2)}_a$ of regularity $\vert a\vert-2$. Keeping in mind that the expression $(4\alpha-2)$ may change from line to line, this yields
	\begin{align*}
	\P_{\varepsilon(u,\cdot)}Lu &=  \P_{\varepsilon(u,\cdot)}\P_{u_a}\xi^{(2)}_a + (4\alpha-2)   \\
										&= \P_{\varepsilon(u,\cdot)u_a}\xi_a^{(2)} + {\sf R}^\circ\big(\varepsilon(u,\cdot), u_a, \xi_a^{(2)}\big) + (4\alpha-2)   \\
										&= \P_{\varepsilon(u,\cdot)u_a}\xi_a^{(2)} + {\sf R}^\circ\big(\varepsilon(u,\cdot), u_\tau, L\tau\big) + (4\alpha-2)   \\
										&= \P_{\varepsilon(u,\cdot)u_a}\xi_a^{(2)} + \DR^\circ\big(\varepsilon(u,\cdot)u_{\tau\Mj},\Mj,L\tau\big) + (4\alpha-2)   \\
										&= \P_{\varepsilon(u,\cdot)u_a}\xi_a^{(2)} + \P_{d_0^{-1}d'(u)u_\Mk u_{\tau\Mj} + \varepsilon(u,\cdot)u_{\tau\Mj\Mk}}\DR^\circ\big(\Mk,\Mj,L\tau\big) + (4\alpha-2)   \\
										&= \P_{\varepsilon(u,\cdot)u_a}\zeta_a^{(2)} + \P_{d_0^{-1}d'(u)u_\Mk u_{\tau\Mj}}Y_{\Mk,\tau\Mj}^d + (4\alpha-2)
	\end{align*}
	The term ${\sf P}_{\varepsilon(u,\cdot)u_{\tau\sigma\gamma}}{\sf R}^\circ(\gamma,\sigma,L\tau)$ has been added to ${\sf P}_{\varepsilon(u,\cdot)u_a}\xi_a^{(2)}$, with $a=\tau\sigma\gamma$, resulting in changing $\xi^{(2)}_a$ to $\zeta^{(2)}_a$. We rewrite this formula under the form
	\begin{align*}
	\P_{\varepsilon(u,\cdot)}Lu = \sum_{\tau\in\mathscr{T}}\P_{\varepsilon(u,\cdot)u_\tau}L\tau + \sum_{a\notin\mathscr{T}} \P_{\varepsilon(u,\cdot)u_a}\zeta_a^{(2)} + \P_{d_0^{-1}d'(u)u_\Mk u_{\tau\Mj}}Y_{\Mk,\tau\Mj}^d + (4\alpha-2),
	\end{align*}
	to put forward the terms $\P_{\varepsilon(u,\cdot)u_\tau}L\tau$, of regularity $\alpha-2$. This is the only term in the right hand side of equation \eqref{EqReformulation} that has the same regularity as the noise $\zeta$.   \vspace{0.15cm}

	\item The terms 
	$$
	\P_{Lu}\varepsilon(u,\cdot) = \P_{Lu}d_0^{-1}d(u) - \P_{Lu}{\bf 1}  = \P_{Lu}\big(d_0^{-1}d(u)\big) + (4\alpha-2)
	$$ 
	and 
	$$
	\PI\big(\varepsilon(u,\cdot), Lu\big) = \PI\big(d_0^{-1}d(u), Lu\big) - \PI({\bf 1},Lu) = \PI\big(d_0^{-1}d(u), Lu\big) + (4\alpha-2)
	$$
	are dealt with using the correctors ${\sf C},{\sf C}_L^<,{\sf C}_L^>$ and ${\sf C}_L$, to take care of paraproducts, and their refined versions $C_{L,(1)}^1$, etc., to take care of remainder terms in paracontrolled expansions. Recall the $\sf E$-type form of the continuity statements on these operators. Recall also that $d_0$ and $d_0^{-1}$ are smooth. Using the $\sf E$-notation for operators of $\sf E$-type, such as in the introduction of Section \ref{SectionPCScene}, we have 
	\begin{align*}
	\P_{Lu}\big(d_0^{-1}d(u)\big) + \PI\big(d_0^{-1}d(u), Lu\big) &= {\sf E}^{-2}\big(d_0^{-1}d(u), u\big)   \\
	&= d_0^{-1}d'(u)E^{-2}(u,u) + d_0^{-1}d^{(2)}(u)E^{-2}(u,u,u) + (4\alpha-2).
	\end{align*}
The analysis of the term ${\sf E}^{-2}(u,u)$ is conveniently done as follows. (This computation was already done at length in \cite{BB3}.) We first write it in multiplicative form
\begin{equation*} \begin{split}
{\sf E}^{-2}(u,u) &= u_{\tau_1}{\sf E}^{-2+\vert\tau_1\vert}(u) + {\sf E}^{-2+\vert\tau_1\vert}(u_{\tau_1}, u) + (5\alpha-2)   \\
				   &= \Big\{u_{\tau_1}u_{\tau_2}{\sf E}^{-2+\vert\tau_1\vert+\vert\tau_2\vert} + u_{\tau_1}{\sf E}^{-2+\vert\tau_1\vert+\vert\tau_2\vert}(u_{\tau_2}) + (5\alpha-2)\Big\}   \\
				   &\quad+ \Big\{u_{\tau_1\sigma_1}{\sf E}^{-2+\vert\tau_1\vert+\vert\sigma_1\vert}(u) +  {\sf E}^{-2+\vert\tau_1\vert+\vert\sigma_1\vert}(u_{\tau_1\sigma_1}, u) + (5\alpha-2) \Big\}   \\
				   &\quad+ (5\alpha-2)   \\
				   &= \Big\{u_{\tau_1}u_{\tau_2}{\sf E}^{-2+\vert\tau_1\vert+\vert\tau_2\vert} + u_{\tau_1}u_{\tau_2\sigma_2}{\sf E}^{-2+\vert\tau_1\vert+\vert\tau_2\vert+\vert\sigma_2\vert}  \\
				   &\qquad+ u_{\tau_1}u_{\tau_2\sigma_2\mu_2}{\sf E}^{-2+\vert\tau_1\vert+\vert\tau_2\vert+\vert\sigma_2\vert+\vert\mu_2\vert} + (5\alpha-2)\Big\}   \\
				   &\quad+ \Big\{u_{\tau_1\sigma_1}{\sf E}^{-2+\vert\tau_1\vert+\vert\sigma_1\vert+\vert\tau_2\vert} + u_{\tau_1\sigma_1}u_{\tau_2\sigma_2}{\sf E}^{-2+\vert\tau_1\vert+\vert\sigma_1\vert+\vert\tau_2\vert+\vert\sigma_2\vert} + (5\alpha-2)   \\
				   &\qquad+ u_{\tau_1\sigma_1\mu_1}u_{\tau_2}{\sf E}^{-2+\vert\tau_1\vert+\vert\sigma_1\vert+\vert\mu_1\vert+\vert\tau_2\vert} + (5\alpha-2) \Big\}   \\
				   &\quad+ (5\alpha-2).
\end{split} \end{equation*}
Each term above that is not a remainder $(5\alpha-2)$ is of the form 
$$
(\star)\,{\sf E}^\beta = {\sf P}_{(\star)}{\sf E}^\beta + {\sf F}^\beta(\star) + {\sf E}^\beta(\star),
$$
for different values of $\beta$, and $(\star)$ either of the form $u_a$ or $u_au_b$, with $a,b\in\mathscr{A}$. The term ${\sf P}_{(\star)}{\sf E}^\beta$ has the expected form. We use the paracontrolled structure of $u_a$ and the $\sf F$-expansion property to deal with ${\sf F}^\beta(u_a)$. To deal with ${\sf F}^\beta(u_au_b)$, write first
$$
{\sf F}^\beta(u_au_b) = {\sf F}^\beta\big({\sf P}_{u_a}u_b\big) + {\sf F}^\beta\big({\sf P}_{u_b}u_a\big) + {\sf F}^\beta\big({\sf \Pi}(u_a,u_b)\big),
$$
and use the $\sf F$-expansion property for the first two terms. For the resonant term, we use the commutator operator $\sf D$ and its continuity properties, recalled in Appendix \ref{AppendixPCOverview}, to expand first the resonant term in the form 
$$
{\sf \Pi}(u_a,u_b) = {\sf P}_{u_{a\tau}}{\sf \Pi}(\tau, u_b) + {\sf D}(u_{a\tau},\tau,u_b),
$$
and then expand the paraproduct inside the operators ${\sf \Pi}$ and $\sf D$, using the paracontrolled forms of $u_b$ and $u_{a\tau}$. We leave the details to the reader; all these operations are only done up to remainders of positive regularity $5\alpha-2$. We also leave the analysis of the term ${\sf E}^{-2}(u,u)$ to the reader. These computations give in the end
	\begin{align*}	
	\P_{Lu}\big(d_0^{-1}d(u)\big) + \PI\big(d_0^{-1}d(u), Lu\big) = \P_{d_0^{-1}d'(u)u_au_b}\zeta_{ab}^{(2)} + \P_{d_0^{-1}d^{(2)}(u)u_au_bu_c}\zeta_{abc}^{(2)} + (4\alpha-2).   \vspace{0.1cm}
	\end{align*}
	
	\item For the terms involving the vector fields $a_i(u,\cdot)V_iu$, we simply note that
	$$
	\P_{V_iu}a_i(u,\cdot)+\PI(a_i(u,\cdot),V_iu)=(2\alpha-1)=(4\alpha-2),
	$$
	since $2\alpha-1>4\alpha-2$, and that 
	\begin{align*}
	V_iu&= V_i\PT_{u_\tau}\tau + (4\alpha-2)   \\
	&= \P_{u_\tau}V_i\tau + \DV_i(u_\tau,\tau) + (4\alpha-2)   \\
	&= \P_{u_\tau}V_i\tau + (4\alpha-2).
	\end{align*}
	using Theorem \ref{ThmVi2}.
\end{itemize}
\end{Dem}

\medskip

We insist again on the fact that all the implicit sums on repeated indices above converge as a consequence of the bound \eqref{EqConvergenceCondition} satisfied by paracontrolled systems, and from the continuity estimates from Section \ref{SectionAdditionalCorrectors}.

\medskip

\textbf{\textsf{Step 3.}} We did not say so far which reference set $\mathscr{T}$ choosing in \textbf{\textsf{Step 1}}. We build $\mathscr{T}$ from the fixed point formulation of equation \eqref{EqReformulation} from Proposition \ref{PropDecomposition}. 

\ssk

One identifies from equation \eqref{EqDecomposition} a number of constraints that $\mathscr{T}$ needs to satisfy. Denote by $e=(a_1,\dots,a_k)$ a generic sentence with words in $\mathscr{A}$, with $\vert e\vert := \vert a_1\vert + \cdots + \vert a_k\vert$.
\begin{equation} \begin{split}   \label{EqConditionScrT1}
\left\{
    \begin{array}{rlll}
	\mathscr{L}^{-1}(\zeta)&\in\mathscr{T}_1, &   \\
	(\mathscr{L}^{-1}L)(\mathscr{T}_i)&\subset \mathscr{T}_i, &\textrm{for } 1\leq i\leq 3,  \\
	\mathscr{L}^{-1}\big(\zeta^{(1)}_e\big)&\subset \mathscr{T}_{i+1}, &\textrm{for } \vert e\vert=i\alpha\leq 2\alpha,   \\
	\mathscr{L}^{-1}\big(\zeta^{(2)}_e\big)&\subset \mathscr{T}_i, &\textrm{for } \vert e\vert=i\alpha\leq 3\alpha,  \textrm{ and } e\notin\mathscr{T}, \\
	\mathscr{L}^{-1}\big(\zeta_j(\mathscr{T}_1)\big)&\subset \mathscr{T}_3. &
    \end{array}
\right.
\end{split} \end{equation}
Recall from Appendix \ref{AppendixPCOverview} the definition of the operator $\sf R$. Requiring 
\begin{equation}  \label{EqConditionScrT2}
{\sf \Pi}(\tau,\sigma)\in\mathscr{T}_2 \textrm{ and } {\sf R}({\bf 1}, \tau,\sigma)\in\mathscr{T}_2, \quad \forall \tau,\sigma\in\mathscr{T}_1.
\end{equation}
ensures that for $u$ paracontrolled to order $3$ by the a reference set $\mathscr{T}$, all the functions $f(u), f'(u)u_a, f^{(2)}(u)u_au_b$, etc. that appear as arguments of the paraproducts in identity \eqref{EqDecomposition}, have a second order paracontrolled expansion with respect to that reference set $\mathscr{T}$. We define $\mathscr{T} = \mathscr{T}_1\cup\mathscr{T}_2\cup\mathscr{T}_3$, as the smallest set of reference functions satisfying the constraints \eqref{EqConditionScrT1} and \eqref{EqConditionScrT2}. This construction recipe for $\mathscr{T}$ gives back the finite set $\mathscr{T}^\circ$ used for the study of the semilinear generalised (PAM) equation in \cite{BB3}, if one replaces the preceding infinite set $\mathscr{T}_1$ be the one point set $\big\{\mathscr{L}^{-1}\zeta\big\}$. In a sense, one can see $\mathscr{T}^\circ$ as the `skeleton' of $\mathscr{T}$, each occurence of $\mathscr{L}^{-1}(\zeta)$ in an element of $\mathscr{T}^\circ$ being possibly any element of $\mathscr{T}_1$ in $\mathscr{T}$. Given $\tau\in\mathscr{T}$, denote by $n_\tau$ the total number of times that the operator $\mathscr{L}^{-1}L$ appears in the formal expression for $\tau$.

\medskip

\noindent \textbf{\textsf{Assumption (A).}} {\it  There exists positive constants $k$ and $C>1$ such that one has
$$
\|\tau\|_{\mcC^{|\tau|}} \leq k\,C^{n_\tau},
$$
for all $\tau\in\mathscr{T}$.   }

\medskip

See remark 2 after the proof of Theorem \ref{ThmMain} for comments on this assumption. With that choice of $\mathscr{T}$, given a system $\widehat{u}$ paracontrolled by $\mathscr{T}$, the function
$$
\mathscr{L}^{-1}\left(f(u)\zeta + \varepsilon(u,\cdot)Lu + \sum_{i=1}^\ell a_i(u,\cdot)V_iu\right)
$$
is the first element of a system paracontrolled by $\mathscr{T}$ that we denote by
$$
\Psi\big((u_a)_{a\in\mathscr{A}}\big).
$$
Write
$$
\Phi\big((u^\sharp_a)_{a\in\mathscr{A}}\big)
$$
for the associated map that gives the collection of all the remainders in the paracontrolled expansion of the different elements of $\Psi\big((u_a)_{a\in\mathscr{A}}\big)$. Note that the fixed point identity 
$$
u = \sum_{\tau\in\mathscr{T}} \widetilde{P}_{u_\tau}\tau + u^\sharp = \mathscr{L}^{-1}\Big(f(u)\zeta + \varepsilon(u,\cdot)Lu + \sum_{i=1}^\ell a_i(u,\cdot)V_iu\Big) + Pu_0
$$
identifies then each $u_\tau$ in the left hand side to an explicit function $h_\tau(u)$ of $u$ only. One has for instance 
$$
u_\tau = \varepsilon(u,\cdot)^kf(u), \textrm{ for } \tau= (\mathscr{L}^{-1}L)^k\big(\mathscr{L}^{-1}\zeta\big),
$$

\ssk

{\it We now choose the exponents $(\beta_a)_{a\in\mathscr{A}}$ in $(5/2,\alpha)$ in such a way that $\beta_a>\beta_{a'}$, if the word $a$ has more letters that $a'$, and $\beta_a>\beta_{a'}$, if $a$ and $a'$ have the same number of letters and $|a|<|a'|$. Given the above skeleton picture of $\mathscr{T}$, this can be done in such a way that the $\beta_a$ take only finitely many values. This will be import in order to prove that the map $\Phi$ is a contraction for $T$ small enough.}

\bigskip

Denote by $\widehat{u}^\sharp = (u_a^\sharp)_{a\in\mathscr{A}}$ a generic element of the product space
$$
\prod_{a\in\mathscr{A}} \mcC^{3\alpha+\beta_a-|a|},
$$ 
endowed with the norm
$$
\vvvert \widehat{u}^\sharp\vvvert := \sum_{a\in\mathscr{A}}\|u_a^\sharp\|_{\mcC^{3\alpha+\beta_a-|a|}}\llparenthesis a\rrparenthesis.
$$
Given $u_0\in C^{4\alpha}$, set $h_\emptyset(u_0) := u_0$, and define
$$
\mcS(u_0) := \left\{\widehat{u}^\sharp;\,\vvvert\widehat{u^\sharp}\vvvert<\infty, \textrm{ and } {u_a^\sharp}_{|_{t=0}}=h_a(u_0),\;\forall\,a\in\mathscr{A}\right\};
$$
this is a closed subspace of $\Big(\prod_{a\in\mathscr{A}} \mcC^{3\alpha+\beta_a-|a|}, \vvvert\cdot\vvvert\Big)$.

\medskip

\begin{thm} \label{ThmMain}
The map $\Phi$ is a contraction of $\mcS(u_0)$, provided the positive time horizon $T$ is small enough.
\end{thm}

\medskip

This statement means that equation \eqref{EqQgPAM} has a unique local in time solution in the space $\mcS(u_0)$; it depends continuously on $\mathscr{T}$.

\medskip

\begin{Dem}
Recall we use exclusively the symbols $\tau,\sigma$ for letters from the alphabet $\mathscr{T}$, while we write $a,b,c$ for elements of $\mathscr{A}$ -- possibly words with only one lettres.

\ssk

$\bullet$ \textit{We first prove that $\Phi$ is a well-defined map from $\mcS(u_0)$ into itself}, that is show is that the condition
$$
\sum_{a\in\SA}\|u_a^\sharp\|_{\CC^{3\alpha-|a|+\beta_a}}\llparenthesis a\rrparenthesis < \infty
$$
is stable by $\Phi$. We decompose this sum according to the value of $\vert a\vert$.

\ssk

-- For $\vert a\vert=3\alpha$, one has $v_a=v_a^\sharp\in\mcC^{\beta_a}$, and the condition reads
$$
\sum_{|a|=3\alpha}\|v_a\|_{\mcC^{\beta_a}}\llparenthesis a\rrparenthesis < \infty.
$$
We read on formula \eqref{EqDecomposition} the different possibilities for $a$, of the form $\mathscr{L}^{-1}(\zeta^{(1)}_e)$, with $e\in\{a,(b,c)\}_{a,b,c\in\mathscr{A}}$ and $\vert a\vert=2\alpha$ or $\vert bc\vert=2\alpha$, etc. If for instance $a=\mathscr{L}^{-1}(\zeta^{(1)}_{a'})$ with $|a'|=2\alpha$, we need to show that
$$
\sum_{|a'|=2\alpha}\big\|g'(u)u_{a'}\big\|_{\mcC^{\beta_{a'}}}\big\|\mathscr{L}^{-1}(\zeta^{(1)}_{a'})\big\|_{\mcC^{2\alpha}} < \infty.
$$
This can be seen from a direct computation
\begin{align*}
\sum_{|a'|=2\alpha} \big\|&f'(u)u_{a'}\big\|_{\mcC^{\beta_{a'}}} \big\|\mathscr{L}^{-1}(\zeta^{(1)}_{a'})\big\|_{\mcC^{2\alpha}}   \\
&= \sum_{|\tau| = 2\alpha} \big\|f'(u)u_\tau\big\|_{\mcC^{\beta_{a'}}}\big\|\mathscr{L}^{-1}(\zeta^{(1)}_\tau)\big\|_{\mcC^{2\alpha}} + \sum_{|\Mj|=|\Mk|=\alpha} \|f'(u)u_{\Mj\Mk}\|_{\mcC^{\beta_{a'}}} \big\|\mathscr{L}^{-1}(\zeta^{(1)}_{\sigma\gamma})\big\|_{\mcC^{2\alpha}}   \\
&\lesssim\|f'(u)\|_{\mcC^{\alpha}} \left(\sum_{|\tau|=2\alpha}\|u_\tau\|_{\mcC^{\beta_{a'}}} \|a'\| + \sum_{|\sigma|=|\gamma|=\alpha}\|u_{\sigma\gamma}\|_{\mcC^{\beta_{a'}}}\|\sigma\gamma\|\right)   \\
&\lesssim \sum_{|\tau|=2\alpha}\|u_\tau\|_{\beta_\tau}\|\tau\| + \sum_{|\sigma|=|\gamma|=\alpha} \|u_{\sigma\gamma}\|_{\mcC^{\beta_{\sigma\gamma}}}\|\sigma\gamma\|   \\
&\lesssim \vvvert\widehat{u}^\sharp\vvvert < \infty,
\end{align*}
using that $\beta_\sigma > \beta_{a'}$ since $|\sigma|<|a'|$ and $\beta_{\sigma\gamma}>\beta_{a'}$. We let the reader check the other cases.

\ssk

-- For $|a|=2\alpha$, we need to show
$$
\sum_{|a|=2\alpha}\|v_a^\sharp\|_{\alpha+\beta_a}\llparenthesis a\rrparenthesis < \infty
$$
so we need to compute the remainders $v_a^\sharp$ for all such $a$, which are given by
$$
v_a = \sum_{|\tau|=\alpha}\PT_{v_{a\tau}}\tau + v_a^\sharp.
$$
Here again, different cases can happen depending on $a$. If for instance $a = \mathscr{L}^{-1}(\zeta^{(1)}_{a'})$, with $|a'|=\alpha$, so $a'=\sigma\in\mathscr{T}_1$, we have $v_a=f'(u)u_\sigma$, and 
\begin{align*}
f'(u)u_\sigma &= \P_{f'(u)}u_\Mj + \P_{u_\Mj}f'(u)+\PI\big(f'(u),u_\Mj\big)   \\
&= \PT_{f'(u)u_{\sigma\gamma} + f^{(2)}(u)u_\sigma u_\gamma}\gamma + \Big\{\DR\big(f'(u),u_{\sigma\gamma},\gamma\big) + \DR\big({\bf 1},f '(u)u_{\sigma\gamma},\gamma\big) + \DR\big(f^{(2)}(u)u_\sigma, u_\gamma,\gamma\big)   \\
&\quad+ \DR({\bf 1}, f^{(2)}(u)u_\sigma u_\gamma,\gamma\big) + \P_{u_\sigma} f'(u)^\sharp + \PI\big(f'(u), u_\sigma\big) \Big\}  \\
&=: \PT_{v_{a\gamma}}\gamma + v_a^\sharp
\end{align*}
where all term in the remainder $v_a^\sharp$ satisfy the convergence condition. As an example, we have
\begin{align*}
\sum_{|\Mj|=\alpha}\sum_{|\Mk|=\alpha} \left\|\DR({\bf 1}, f'(u)u_{\Mj\Mk}, \Mk)\right\|_{\alpha+\beta_a}\|\Mj\|&\lesssim\|g'(u)\|_\alpha\sum_{|\Mj|=|\Mk|=\alpha}\|u_{\Mj\Mk}\|_{\beta_{\Mj\Mk}}\|\Mj\Mk\| \lesssim \vvvert\widehat{u}^\sharp\vvvert < \infty.
\end{align*}
The reader is invited to check the other cases.

\ssk

-- A direct computation also shows that
$$
\sum_{|a|=\alpha} \|v_a^\sharp\|_{2\alpha+\beta_a}\llparenthesis a\rrparenthesis < \infty.
$$
The remaining details are left to the reader.

\bigskip

$\bullet$ \textit{We now prove that $\Phi$ is a contraction for a positive time horizon $T$ small enough.} Pick $\widehat{u}^\sharp = (u_a^\sharp)_{a\in\SA}\in\mcS(u_0)$ and $\widehat{v}^\sharp = (v_a^\sharp)_{a\in\SA}\in\mcS(u_0)$. Since both paracontrolled systems are in the solution space,
$$
\widehat{w}:=\Phi(\widehat{u})-\Phi(\widehat{v})
$$
is a system paracontrolled by $\mathscr{T}$ at order $3$, whose elements are equal to $0$ at time $0$. This allows us to gain a factor $T^{(\gamma'-\gamma)/2}$ when comparing the norms of such functions in two different parabolic H\"older spaces with respective exponents $\gamma$ and $\gamma'$. From Propostion \ref{PropDecomposition}, we have
$$
\widehat{w}_\emptyset = \sum_{|\Mi|\le3\alpha}\PT_{w_\Mi}\Mi+w^\sharp
$$
with explicit formulas for $w_\Mi=\Phi(\widehat{u})_\Mi-\Phi(\widehat{v})_\Mi$. In the expansion of $\widehat{w}_\emptyset$, we have to control two types of terms,
that is \textbf{\textsf{(i)}} $\|w_\Mi\|_{\beta_\Mi}$ for $|\Mi|=3\alpha$ and \textbf{\textsf{(ii)}}  $\|w^\sharp\|_{3\alpha+\beta_\emptyset}$.

\ssk

\textbf{\textsf{(i)}} We first consider the terms $w_\Mi$. For example, we have to control
$$
\left\|f'(u)u_a-f'(v)v_a\right\|_{\beta_{\tau}}, \quad\text{ with } |a|=2\alpha,\text{ and }\tau = \mathscr{L}^{-1}(\zeta_a^{(1)}).
$$
and this is done writing
\begin{align*}
\big\|f'(u)u_a - f'(v)v_a\big\|_{\beta_{\Mj}} &\lesssim \big\|\left(f'(u)-f'(v)\right)u_a\big\|_{\beta_\Mj} + \big\|f'(v)(u_a-v_a)\big\|_{\beta_\Mj}   \\
&\lesssim T^{\frac{\alpha-\beta_{\Mj}}{2}}\big\|f'(u) - f'(v)\big\|_\alpha \|u_a\|_{\beta_\Mj} + T^{\frac{\beta_a-\beta_\Mj}{2}} \big\|f'(v)\big\|_\alpha \|u_a-v_a\|_{\beta_a}   \\
&\lesssim T^{\frac{\alpha-\beta_{\Mj}}{2}}\left(\|f\|_{C_b^2}\big(1+\|u\|_\alpha\big)\|u_a\|_{\beta_\Mj}\right) \|u-v\|_\alpha + T^{\frac{\beta_a-\beta_\Mj}{2}}\|f'(v)\|_\alpha \|u_a - v_a\|_{\beta_a}   \\
&\lesssim\left\{T^{\frac{\alpha-\beta_{\Mj}}{2}}\Big(\|f\|_{C_b^2}\big(1+\|u\|_\alpha\big)\|u_a\|_{\beta_\Mj}\Big) + T^{\frac{\beta_a-\beta_\Mj}{2}}\|f'(v)\|_\alpha\right\} \vvvert \widehat{u}^\sharp - \widehat{v}^\sharp\vvvert
\end{align*}
Another example is
$$
\big\|f^{(2)}(u)u_au_b - f^{(2)}(v)v_av_b\big\|_{\beta_\Mj}
$$
with $|a|+|b|=2\alpha$ and $\Mj\in\mathscr{T}$ given by $\SL^{-1}\left(\zeta_{ab}^{(1)}\right)$. It is dealt with writing
\begin{align*}
\big\|f^{(2)}(u)u_au_b - f^{(2)}(v)v_av_b\big\|_{\beta_\Mj} \lesssim \left\|\left(f^{(2)}(u)-f^{(2)}(v)\right)u_au_b\right\|_{\beta_\Mj} + \big\|f^{(2)}(v)(u_au_b-v_av_b)\big\|_{\beta_\Mj}  
\end{align*}
\begin{align*}
&\lesssim T^{\frac{\alpha-\beta_\Mj}{2}}\|f^{(2)}(u)-f^{(2)}(v)\|_\alpha\|u_au_b\|_{\beta_\Mj}+T^{\frac{\min(\beta_a,\beta_b)-\beta_\Mj}{2}} \big\|f^{(2)}(v)\big\|_{\beta_\Mj}\big\|u_au_b - v_av_b\big\|_{\min(\beta_a,\beta_b)}   \\
&\lesssim\left(T^{\frac{\alpha-\beta_\Mj}{2}}\|f\|_{C_b^3}\|u_au_b\|_{\beta_\Mj}+T^{\frac{\min(\beta_a,\beta_b)-\beta_\Mj}{2}}\|f^{(2)}(v)\|_{\beta_\Mj}\right)\vvvert\widehat{u}^\sharp-\widehat{v}^\sharp\vvvert.
\end{align*}
All the other terms are dealed with using the following four inequalities. 
\begin{itemize}
   \item One has
$$
\Big\|\varepsilon(u,\cdot)u_a - \varepsilon(v,\cdot)v_a\Big\|_{\beta_\Mj} \lesssim T^{\frac{\alpha-\beta_\Mj}{2}} \big\|\varepsilon(u,\cdot)-\varepsilon(v,\cdot)\big\|_\alpha \|u_a\|_{\beta_\Mj} + T^{\frac{\beta_a-\beta_\Mi}{2}} \|\varepsilon(v,\cdot)\|_\alpha\|u_a-v_a\|_{\beta_a}
$$
for $a\notin\mathscr{T},|a|=3\alpha$ and $\Mj\in\mathscr{T}$ given by $\SL^{-1}\left(\zeta_a^{(2)}\right)$.   \vspace{0.1cm}
   
   \item One has
   \begin{align*}
\big\|d_0^{-1}d'(u)u_au_b-d_0^{-1}d'(v)&v_av_b\big\|_{\beta_\Mj} \lesssim T^{\frac{\alpha-\beta_\Mj}{2}}\big\| d_0^{-1}d'(u) - d_0^{-1}d'(v)\big\|_\alpha\|u_au_b\|_{\beta_\Mj}   \\
&+ T^{\frac{\min(\beta_a,\beta_b)-\beta\Mj}{2}} \big\|d_0^{-1}d'(v)\big\|_{\beta_\Mj}\big\|u_au_b - v_av_b\big\|_{\min(\beta_a,\beta_b)}
\end{align*}
for $|a|+|b|=3\alpha$ and $\Mj\in\mathscr{T}$ given by $\SL^{-1}\left(\zeta_{ab}^{(2)}\right)$.   \vspace{0.1cm}
   
   \item One has
   \begin{align*}
\big\|d_0^{-1}d^{(2)}(u)u_au_bu_c-d_0^{-1}d^{(2)}(v)&v_av_bv_c\big\|_{\beta_\Mj} \lesssim T^{\frac{\alpha-\beta_\Mj}{2}} \big\|d_0^{-1}d^{(2)}(u)-d_0^{-1}d^{(2)}(v)\big\|_\alpha \|u_au_bu_c\|_{\beta_\Mj}   \\
&+ T^{\frac{\min(\beta_a,\beta_b,\beta_c)-\beta_\Mi}{2}} \big\|d_0^{-1}d^{(2)}(v)\big\|_{\beta_\Mj} \big\|u_au_bu_c-v_av_bv_c\big\|_{\min(\beta_a,\beta_b,\beta_c)}
\end{align*}
for $|a|+|b|+|c|=3\alpha$ and $\Mj\in\mathscr{T}$ given by $\SL^{-1}\left(\zeta_{abc}^{(2)}\right)$.   \vspace{0.1cm}
   
   \item One has
   $$
\big\| d_\ell(u,\cdot)u_\Mi - d_\ell(v,\cdot)v_\Mj\big\|_{\beta_\Mj} \lesssim T^{\frac{\alpha-\beta_\Mj}{2}} \big\| d_\ell(u,\cdot) - d_\ell(v,\cdot)\big\|_\alpha\|u_\Mi\|_{\beta_\Mj} + T^{\frac{\beta_\Mi-\beta_\Mj}{2}} \big\| d_\ell(v,\cdot)\big\|_{\beta_\Mj} \big\| u_\Mi - v_\Mi\big\|_{\beta_\Mi}
$$
for $|\Mi|=\alpha$ and $\sigma\in\mathscr{T}$ given by $\mathscr{L}^{-1}(\zeta_{j,\tau})$.
\end{itemize}

\medskip

There is only one case where we use that $\varepsilon(u,\cdot)$ is small when $T$ is small, to estimate
\begin{equation} \label{EqEstimateVarepsilon}
\Big\|\varepsilon(u,\cdot)u_\Mi - \varepsilon(v,\cdot)v_\Mi\Big\|_{\beta_\Mj} \lesssim \big\|\left(\varepsilon(u,\cdot)-\varepsilon(v,\cdot)\right)u_\Mi\big\|_{\beta_\Mj} + \big\|\varepsilon(v,\cdot)(u_\Mi-v_\Mi)\big\|_{\beta_\Mj}.
\end{equation}
While we have
$$
\|\left(\varepsilon(u,\cdot)-\varepsilon(v,\cdot)\right)u_\Mi\|_{\beta_\Mj}\lesssim T^{\frac{\alpha-\beta_\Mj}{2}} \big\|\varepsilon(u,\cdot) - \varepsilon(v,\cdot)•big\|_\alpha \|u_\Mi\|_{\beta\Mj}
$$
we do not gain a $T$-dependent fact using the regularity of $u_\Mi-u_\Mj$ in the second term of the right hand side of inequality \eqref{EqEstimateVarepsilon}, since $\beta_\Mi=\beta_\Mj$. We use instead that
\begin{align*}
\big\|\varepsilon(v,\cdot)(u_\Mi-v_\Mi)\big\|_{\beta_\Mj} &\lesssim \|\varepsilon(v,\cdot)\|_{\beta_\Mi}\|u_\Mi-v_\Mi\|_{\beta_\Mi}   \\
&\lesssim \left(\big\|d(v)-d(u_0)\big\|_{\beta_\Mi} + \big\|d(u_0)-d(\overline{u_0})\big\|_{\beta_\Mi}\right) \|u_\Mi-v_\Mi\|_{\beta_\Mi}   \\
&\lesssim \left(T^{\frac{\beta_\Mi-\alpha}{2}}\big\|d(v)-d(u_0)\big\|_{\alpha} + \big\|d(u_0)-d(\overline{u_0})\big\|_{\beta_\Mi}\right) \|u_\Mi-v_\Mi\|_{\beta_\Mi}
\end{align*}
since $d(v)-d(u_0)$ is equal to $0$ at $t=0$, and the factor $\|d(\overline{u_0})-d(u_0)\|_{\beta_\Mi}$ is as small as we want for $\overline{u_0}$ close enough to $u_0$ in $C^\alpha$.

\ssk

\textbf{\textsf{(ii)}} For the remainder terms, the same arguments with an explicit computations of $(\sharp)(\widehat{w})$ is enough. For example, we have for $|\Mi|=|\Mj|=|\Mj|=\alpha$ a term
$$
\big\|\DL(((u_{\Mi\Mj\Mk},\Mk),\Mj),\Mi)\big\|_{3\alpha+\beta_\emptyset} \lesssim T^{\frac{\beta_{\Mi\Mj\Mk}-\beta_\emptyset}{2}} \| u_{\Mi\Mj\Mk}\|_{\beta_{\Mi\Mj\Mk}}\|\Mi\|_\alpha\|\Mj\|_\alpha\|\Mk\|_\alpha
$$
or
\begin{align*}
\Big\|\PT_{\varepsilon(u,\cdot)}(\SL^{-1}L)u^\sharp &- \PT_{\varepsilon(v,\cdot)}(\SL^{-1}L)v^\sharp\Big\|_{3\alpha+\beta_\emptyset}\lesssim T^{\frac{\alpha-\beta_\emptyset}{2}}\|\varepsilon(u,\cdot)-\varepsilon(v,\cdot)\|_\alpha\|u^\sharp\|_{3\alpha+\beta_\emptyset}   \\
&+ T^{\frac{\alpha-\beta_\emptyset}{2}}\|\varepsilon(v,\cdot)\|_\alpha\|u^\sharp-v^\sharp\|_{3\alpha+\beta_\emptyset}.
\end{align*}

\medskip

The above considerations deal with $\widehat{w}_\emptyset$. The analysis of the general terms $\widehat{w}_a$, with $a\in\mathscr{A}$, is similar or easier; it is left to the reader.

\medskip

In order to get a contraction, we see that the different terms $T^{\frac{\beta}{2}}$ can be chosen smaller than one for $T$ small enough. Since there is a finite numbers of $\beta$'s, this shows that $\Phi$ is indeed a contraction for $T$ small enough, and $\overline{u_0}$ close enough to $u_0$ in $C^\alpha$.
\end{Dem}

\medskip

\noindent \textbf{\textsf{Remarks.}} \textbf{\textsf{1.}} \textit{Condition \eqref{EqConditionScrT1} makes plain sense if the noise $\zeta$ is smooth enough. In a low regularity setting where the noise is a random distribution one needs to build the elements of the space $\mathscr{T}$ by probabilistic means from a renormalisation algorithm. A constraint like $(\mathscr{L}^{-1}L)(\zeta^{(1)}_a)\subset \mathscr{T}_{i+1}$, for $\vert a\vert=i\alpha$, is then understood as the requirement that the inclusion holds for each regularisation parameter, for the corresponding function/distribution spaces, and each regularized/renormalised function/distribution is converging in a probabilistic sense to a limit in its natural space.}   \vspace{0.15cm}

\textbf{\textsf{2.}} \textit{One expects that a version of the \emph{BPHZ} renormalisation works for the renormalisation of semilinear singular stochastic PDEs. (While this has been fully proved in the regularity structures setting \cite{BHZ, CH}, this remains to be worked out in a paracontrolled setting. See \cite{BH1} for a strong hint that this is indeed the case.) On that basis, it is most likely that an ad hoc renormalisation process for the quasilinear setting would be given by the same renormalisation process as in the semilinear setting, with trees with branches of ad hoc length used instead of their skeleton `semilinear' trees. As a consequence, one expects an estimate of the form
$$
\|\tau\|_{\mcC^{|\tau|}} \leq k_{\tau^\circ}\,C^{n_\tau},
$$
with $\tau^\circ$ the skeleton tree corresponding to $\tau$ in the semilinear setting and $C>1$ a constant depending only on the operator $\mathscr{L}^{-1}L$. As there are only finitely many trees $\tau^\circ$ in a subcritical regime, one could take a uniform constant $k$ instead of $k_{\tau^\circ}$. This point will be the object of a forthcoming work}.   \vspace{0.15cm}

\textbf{\textsf{3.}} \textit{The quasilinear \emph{(gPAM)} equation dealt with in \cite{BDH} involved a {\sl space} white noise $\zeta$ on the two-dimensional torus. One can use in that setting a one point set 
$$
\mathscr{T}_1 = \big\{L^{-1}(\zeta)\big\} =: \{Z_1\},
$$
in the paracontrolled structure, as one has
$$
\mathscr{L}Z_1 = \zeta = LZ_1.
$$
(Recall the computations from the Introduction section.) Since $\mathscr{T}$ reduces to $\mathscr{T}_1$ in the two-dimensional setting, this simplifies the analysis. Even if working with a space white noise, one would need an infinite set $\mathscr{T}_2$ on the $3$-dimensional torus, as one cannot choose only time-independent elements in $\mathscr{T}_2$, so a single time-dependent element in $\mathscr{T}_2$ produces an infinity of them, from the stability condition \eqref{EqConditionScrT1}.}

\bigskip

\appendix
\section{Basics on high order paracontrolled calculus}
\label{AppendixPCOverview}

We recall in this appendix a number of results from \cite{BB2, BB3} that we use in this work. This should help the reader understanding the computations of Appendix \ref{AppendixContinuity} and their mechanics.

We first describe the approximation operators where the heat semigroup plays the role of Fourier theory and Paley-Littlewood projectors in our general geometric setting. The parabolic Hölder space are defined from these operators. We also recall the form of our space-time paraproducts and give some of the continuity estimates on different correctors/commutators and their iterated versions.

\medskip

Recall that we denote by $M$ a $3$-dimensional closed Riemannian manifold and set $\mathcal{M} := [0,T]\times M$, for a finite positive time horizon $T$.  We denote by $e=(\tau,x)$ a generic spacetime point. Denote by $\mu$ the Riemannian volume measure and define the parabolic measure
$$
\nu := dt\otimes\mu.
$$

\bigskip

\subsection{Approximation operators and parabolic Hölder spaces}

In the flat setting of the torus, we can use Fourier theory to approximate Schwartz distributions by smooth functions. We have
$$
f=\lim_{n\to\infty}S_n(f)=\sum_{i\ge-1}\Delta_i(f)
$$
with $\Delta_j$ the Paley-Littlewood projectors. Refer e.g. to \cite{BCD} for basics on Littlewood-Paley theory. Using the heat semigroup, one has in a more general geometric framework
$$
f=\lim_{t\to0}P_t^{(b)}f=\int_0^1 Q_t^{(b)}f\frac{\drm t}{t}+P_1^{(b)}f
$$
where 
$$
Q_t^{(b)}:=\frac{(tL)^be^{-tL}}{(b-1)!}\quad\text{and}\quad-t\partial_tP_t^{(b)}:=Q_t^{(b)}
$$
with $P_0=\text{Id}$. One can show that there exists a polynomial $p_b$ of degree $(b-1)$ such that $P_t^{(b)}=p_b(tL)e^{-tL}$ and $p_b(0)=1$. The operators $Q_t^{(b)}$ and $P_t^{(b)}$ play the role of Paley-Littlewood projector and Fourier series, respectively. Indeed, if one works on the torus, then
$$
\widehat{Q_t^{(b)}}(\lambda)=\frac{\left(t|\lambda|^2\right)^b}{(b-1)!}e^{-|\lambda|^2t}\quad\text{and}\quad\widehat{P_t^{(b)}}(\lambda)=p_b\left(t|\lambda|^2\right)e^{-|\lambda|^2t}
$$
so we see that $Q_t^{(b)}$ localize in frequency around the annulus $|\lambda|\sim t^{-\frac{1}{2}}$ and $P_t^{(b)}$ localize in frequency on the ball $|\lambda|\lesssim t^{-\frac{1}{2}}$. Since the measure $dt/t$ gives unit mass to each interval $[2^{-(i+1)},2^{-i}]$, the operator $Q_t^{(b)}$ is a multiplier that is approximately localized at `frequencies' of size $t^{-\frac{1}{2}}$. However, this decomposition using a continuous parameter does not satisfy the perfect cancellation property $\Delta_i\Delta_j=0$ for $|i-j|>1$, but the identity
$$
Q_t^{(b)}Q_s^{(b)}=\left(\frac{ts}{(t+s)^2}\right)^bQ_{t+s}^{(2b)}
$$
for any $s,t\in(0,1)$. The parameter $b$ encodes a `degree' of cancellation. In order to deal with time approximation, define for $m\in L^1(\IR)$, with support in $\mathbb{R}_+$, the convolution operator
$$
m^\star(f)(\tau):=\int_0^\infty m(\tau-\sigma)f(\sigma)\drm\sigma\quad\text{and}\quad m_t(\cdot):=\frac{1}{t}m\left(\frac{\cdot}{t}\right)
$$
for $\tau\in\IR$ and a positive scaling parameter $t$. Given $I=(i_1,\ldots,i_n)\in\{1,\ldots,\ell\}^n$, define the $n^\textrm{th}$-oder differential operator
$$
V_I:=V_{i_n}\ldots V_{i_1}.
$$
We say that a familly $(\CQ_t)_{t\in(0,1]}$ is Gaussian if each the kernel of each $\CQ_t$ is bounded pointwisely by $\CG_t$, the reference Gaussian kernel. We do not recall its expression here and refer the reader to Section 3.2 of \cite{BB2}.

\medskip

\begin{defn*}
Let $a\in\llbracket0,2b\rrbracket$. We define the \textbf{\textsf{standard collection $\mathsf{StGC}^a$ of operators with cancellation of order}} $a$ as the set of families
$$
\left(\big(t^{\frac{|I|}{2}}V_I\big)\big(tL\big)^{\frac{j}{2}}P_t^{(c)}\otimes\varphi_t^\star\right)_{t\in(0,1]}
$$
where $a=|I|+j+2k$, $c\in\llbracket1,b\rrbracket$ and $\varphi$ a smooth function supported in $[2^{-1},2]$ with bounded first derivative by $1$ such that 
$$
\int \tau^im(\tau)\drm \tau=0\quad\text{for every }0\le i\le k-1.
$$
These operators are uniformly bounded in $L^p(\mcM)$ for every $p\in[1,\infty]$, as functions of the parameter $t\in(0,1]$. We also set
$$
\mathsf{StGC}^{[0,2b]}:=\bigcup_{0\le a\le 2b}\mathsf{StGC}^a.
$$
\end{defn*}

\medskip

A standard family of operator $\CQ\in\mathsf{StGC}^a$ can be seen as a bounded map $t\mapsto\CQ_t$ from $(0,1]$ to the space of bounded linear operator on $L^p(\mcM)$. Since $V_iV_j\neq V_jV_i$, the operators $V_i$ do not commute with $L$ so that
$$
V_IL^be^{-tL}\neq L^be^{-tL}V_I.
$$
For the next proposition, we introduce the notation
$$
\Big(V_I\psi(L)\Big)^\bullet:=\psi(L)V_I
$$
for any holomorphic function $\psi$. This notation is not related to any notion of duality.

\medskip

\begin{prop}   \label{cancel}
Consider $\CQ^1\in\mathsf{StGC}^{a_1}$ and $\CQ^2\in\mathsf{StGC}^{a_2}$ two standard collections with cancellation. Then for every $s,t\in(0,1]$, the composition $\CQ_s^1\circ\CQ_t^{2\bullet}$ has a kernel pointwisely bounded by
\begin{align*}
\Big|K_{\CQ_s^1\circ\CQ_t^{2\bullet}}(e,e')\Big| &\lesssim \left\{\left(\frac{s}{t}\right)^{\frac{a_1}{2}}{\bf 1}_{s<t} + \left(\frac{t}{s}\right)^{\frac{a_2}{2}}{\bf 1}_{s\ge t}\right\}\CG_{t+s}(e,e')   \\
&\lesssim \left(\frac{ts}{(s+t)^2}\right)^{\frac{a}{2}}\CG_{t+s}(e,e')
\end{align*}
with $a:=\min(a_1,a_2)$.
\end{prop} 

\medskip

We also need operators that are not in the standard form but still have some cancellation property.

\medskip

\begin{defn*}
Let $a\in\llbracket0,2b\rrbracket$. We define the collection $\mathsf{GC^a}$ of \textbf{\textsf{operators with cancellation of order $a$}} as the set of families of Gaussian operators $\CQ$ such as the following property holds. For every $s,t\in(0,1]$ and every $\sf S\in\mathsf{StGC}^{a'}$ with $a<a'\le 2b$, the composition $\CQ_s\circ\sf S_t^\bullet$ has a kernel pointwisely bounded by
$$
\left|K_{\CQ_s\circ\sf S_t^\bullet}(e,e')\right|\lesssim\left(\frac{ts}{(t+s)^2}\right)^{\frac{a}{2}}\CG_{t+s}(e,e').
$$
\end{defn*}

\medskip

\begin{defn*}
Given any $\alpha\in(-3,3)$, we define the \textbf{\textsf{parabolic Hölder spaces $\CC^\alpha(\mathcal{M})$}} as the set of distribution $f\in\CD'(\mathcal{M})$ such that
$$
\|f\|_{\CC^\alpha}:=\left\|e^{-L}f\right\|_{L^\infty}+\sup_{\underset{|\alpha|<k\le 2b}{\mcQ\in\mathsf{StGC}^k}}\sup_{t\in(0,1]}t^{-\frac{\alpha}{2}}\|\mcQ_tf\|_{L^\infty}<\infty.
$$
\end{defn*}

\bigskip

\subsection{Parabolic paraproducts, correctors and commutators}

The Paley-Littlewood decomposition can be used to describe a product as
\begin{align*}
fg&=\lim_{n\to\infty}S_n(f)S_n(g)   \\
&=\sum_{i<j-2}\Delta_i(f)\Delta_j(g)+\sum_{|i-j|\le1}\Delta_i(f)\Delta_j(g)+\sum_{i>j+1}\Delta_i(f)\Delta_j(g)   \\
&=\sum_i\Delta_{<i}(f)\Delta_i(g)+\sum_{|i-j|\le1}\Delta_i(f)\Delta_j(g)+\sum_i\Delta_i(f)\Delta_{<i}(g)   \\
&=P_f^0g + \Pi^0(f,g) + P_g^0f,
\end{align*}
the paraproducts $P_f^0g$ and $P_g^0f$ being always well-defined, unlike the resonant term $\Pi^0(f,g)$. In our framework, we use a slightly different identity
\begin{equation} \label{EqDecompositionFG} \begin{split}
fg&=\lim_{t\to0}\CP_t^{(b)}\Big(\CP_t^{(b)}f\cdot\CP_t^{(b)}g\Big)   \\
&=\int_0^1\left\{\CQ_t^{(b)}\big(\CP_t^{(b)}f\cdot\CP_t^{(b)}g\big)+\CP_t^{(b)}\big(\CQ_t^{(b)}f\cdot\CP_t^{(b)}g\big)+\CP_t^{(b)}\big(\CP_t^{(b)}f\cdot\CQ_t^{(b)}g\big)\right\} \frac{dt}{t}   \\
&\hspace{1.5cm}+\CP_1^{(b)}\Big(\CP_1^{(b)}f\cdot\CP_1^{(b)}g\Big).
\end{split} \end{equation}
This corresponds to writing 
$$
fg=\lim_{n\to\infty}S_n\big(S_n(f)S_n(g)\big).
$$
Since $\CP_t^{(b)}$ plays the role of $\Delta_{<i}$ and $\CQ_t^{(b)}$ the role of $\Delta_i$, we want to manipulate this expression to get terms of the following forms
$$
\int_0^1\CP_t^{1\bullet}\left(\CQ_t^1f\cdot\CQ_t^2g\right)\frac{dt}{t},\quad\textrm{or}\quad\int_0^1\CQ_t^{1\bullet}\left(\CQ_t^2f\cdot\CP_t^1g\right)\frac{dt}{t},\quad\text{and}\quad\int_0^1\CQ_t^{1\bullet}\left(\CP_t^1f\cdot\CQ_t^2g\right)\frac{dt}{t},
$$
where $\CQ_1,\CQ_2\in\mathsf{StGC}^c$ encode some cancellation so $c>0$ and $\CP_1\in\mathsf{StGC}^{[0,d]}$ can encode no cancellation. This is done using repeatedly the Leibnitz rule $V_i(fg)=V_i(f)g+fV_i(g)$. For example, we have

{\small \begin{align*}
\int_0^1&\CP_t^{(b)}\left(b^{-1}(tL)\CQ_t^{(b-1)}f\cdot\CP_t^{(b)}g\right)\frac{dt}{t} = b^{-1}\int_0^1\CP_t^{(b)}(tL)\left(\CQ_t^{(b-1)}f\cdot\CP_t^{(b)}g\right)\frac{dt}{t}   \\
&\quad- b^{-1}\int_0^1\CP_t^{(b)}\left(\CQ_t^{(b-1)}f\cdot(tL)\CP_t^{(b)}g\right)\frac{dt}{t} - 2b^{-1}\sum_{i=1}^\ell\int_0^1\CP_t^{(b)}(\sqrt{t}V_i)\left(\CQ_t^{(b-1)}f\cdot(\sqrt{t}V_i)\CP_t^{(b)}g\right)\frac{dt}{t}
\end{align*}   }
where we `take' some cancellation from $\CQ_t^{(b)}$ to the other terms. Starting from identity \eqref{EqDecompositionFG}, repeated use of this decomposition allows to rewrite the product $fg$ as
$$
fg = \P_fg + \PI(f,g) + \P_gf,
$$
where $\P_fg$ is a linear combination of terms of the form 
$$
\int_0^1\CQ_t^{1\bullet}\left(\CP_t^1f\cdot\CQ_t^2g\right)\frac{dt}{t},
$$
and $\PI(f,g)$ is a linear combination of terms of the form
$$
\int_0^1\CP_t^{1\bullet}\left(\CQ_t^1f\cdot\CQ_t^2g\right)\frac{dt}{t},
$$
with $\CQ^1,\CQ^2\in\mathsf{StGC}^{\frac{b}{2}}$ and $\CP^1\in\mathsf{StGC}^{[0,2b]}$, up to the smooth term $\CP_1^{(b)}\left(\CP_1^{(b)}f\cdot\CP_1^{(b)}g\right)$. All the details on this construction and the classical estimates on $\P$ and $\PI$ can be found in Section 4 of \cite{BB2}. In order to solve semilinear singular PDEs, we want to investigate the property of the intertwined operator $\PT$ defined through
$$
\PT_fg=\left(\SL^{-1}\circ\P_f\circ\SL\right)g
$$
for any functions/distributions $f$ and $g$. One can show that $\PT_fg$ is given as a linear combination of
$$
\int_0^1\widetilde\CQ_t^{1\bullet}\left(\CP_t^1f\cdot\CQ_t^2\right)\frac{dt}{t}
$$
with $\widetilde\CQ^1\in\mathsf{GC}^{\frac{b}{4}-2}$, $\CQ^2\in\mathsf{StGC}^{\frac{b}{2}}$ and $\CP^1\in\mathsf{StGC}^{[0,2b]}$. The only difference is that $\widetilde\CQ^1$ is not given by a standard form but still encodes some cancellation. This is however sufficient for $\widetilde{\sf P}$ to enjoy the same continuity properties as $\sf P$.

\medskip 

The study of semilinear singular SPDEs using paracontrolled calculus relies on a number of continuity estimate for different operators, we recall three of them here and refer the reader to \cite{BB3}. Define the $\sf E$-type operators
$$
\DC(a,b,c):=\PI\Big(\PT_ab,c\Big)-a\PI\big(b,c\Big)
$$
and its iterate
$$
\DC\Big((a,b),c,d\Big):=\DC\Big(\PT_ab,c,d\Big)-a\DC\Big(b,c,d\Big).
$$

\medskip

\begin{prop}
\begin{itemize}
	\item Let $\alpha\in(0,1)$ and $\beta,\gamma\in(-3,3)$ such that
	$$
	\beta+\gamma<0\quad\text{and}\quad0<\alpha+\beta+\gamma<1.
	$$
	Then the corrector $\DC$ has a unique extension as a continuous operator from $\CC^\alpha\times\CC^\beta\times\CC^\gamma$ to $\CC^{\alpha+\beta+\gamma}$.   \vspace{0.15cm}
	
	\item Let $\alpha_1,\alpha_2\in(0,1)$ and $\beta,\gamma\in(-3,3)$ such that
	$$
	\alpha_1+\beta+\gamma<0,\quad\alpha_2+\beta+\gamma<0\quad\text{and}\quad0<\alpha_1+\alpha_2+\beta+\gamma<1.
	$$
	Then the iterated corrector $\DC$ has a unique extension as a continuous operator from $\CC^{\alpha_1}\times\CC^{\alpha_2}\times\CC^\beta\times\CC^\gamma$ to $\CC^{\alpha_1+\alpha_2+\beta+\gamma}$.
\end{itemize}
\end{prop}

\medskip

Note that the Hölder regularity exponent of the first argument in the corrector $\DC$ has to be less than $1$ in the above statement. In order to gain more information from a regularity exponent in the interval $(1,2)$, one needs to consider the refined corrector given for any $e\in\mathcal{M}$, by
$$
\DC_{(1)}\Big(a,b,c\Big)(e):=\DC\Big(a,b,c\Big)(e)-\sum_{i=1}^\ell\gamma_i\big(V_ia\big)(e)\PI\Big(\PT_{\delta_i(e,\cdot)}b,c\Big)(e)
$$
where $\delta_i$ is given by 
$$
\delta_i(e,e'):=\chi\big(d(x,x')\big)\langle V_i(x),\pi_{x,x'}\rangle_{T_xM}
$$
with $\chi$ a smooth non-negative function on $[0,+\infty)$ equal to $1$ in a neighbourhood of $0$ with $\chi(r)=0$ for $r\ge r_m$ the injectivity radius of the compact manifold $M$ and $\pi_{x,x'}$ a tangent vector of $T_xM$ of length $d(x,y)$, whose associated geodesic reaches $y$ at time $1$. The functions $\gamma_i$ are defined from the identity
$$
\nabla f = \sum_{i=1}^\ell \gamma_i(V_if)V_i,
$$
for all smoothe real-valued functions $f$ on $M$.

\medskip

\begin{prop}
Let $\alpha\in(1,2)$ and $\beta,\gamma\in(-3,3)$ such that
$$
\alpha+\beta+\gamma>0\quad\text{and}\quad\beta+\gamma<0.
$$
Then the operator $\DC_{(1)}$ has a unique extension as a continuous operator from $\CC^\alpha\times\CC^\beta\times\CC^\gamma$ to $\CC^{\alpha+\beta+\gamma}$.
\end{prop}

\medskip

In the high order setting, we also work with well-defined expressions that are not in the algebraic form  required to set up the fixed point in a simple way; these are the $\sf F$-type terms. We need continuity estimates on
\begin{align*}
\DD(a,b,c) &:= \PI\Big(\PT_ab,c\Big) - \P_a\PI\Big(b,c\Big),   \\
\DR(a,b,c) &:= \P_a\PT_bc - \P_{ab}c,   \\
\DR^\circ(a,b,c) &:= \P_a\P_bc - \P_{ab}c.
\end{align*}

\medskip

\begin{prop}
\begin{itemize}
	\item  Let $\alpha,\beta,\gamma\in(0,3)$. Then the commutator $\DD$ is continuous from $\CC^\alpha\times\CC^\beta\times\CC^\gamma$ to $\CC^{\alpha+\beta+\gamma}$.   \vspace{0.15cm}
	
	\item Let $\beta\in(0,1)$ and $\gamma\in(-3,3)$ such that $\beta+\gamma\in(-3,3)$. Then the operators $\DR$ and $\DR^\circ$ are continuous from $L^\infty\times\CC^\beta\times\CC^\gamma$ to $\CC^{\beta+\gamma}$.   \vspace{0.15cm}
	
	\item Let $\alpha,\beta\in(0,1/2)$ and $\gamma\in(-3,3)$. Then the operator $\DR^\circ$ is continuous from $\CC^\alpha\times\CC^\beta\times\CC^\gamma$ to $\CC^{\alpha+\beta+\gamma}$.
\end{itemize}
\end{prop}

\medskip

We also need continuity estimates on iterates of the operator $\DR^\circ$. However in this case the expansion rule is different depending on which argument we expand.

\medskip

\begin{prop}
\begin{itemize}
	\item Let $\alpha_1,\alpha_2\in(0,1)$ and $\gamma\in(-3,3)$. Then the operator
	$$
	\DR^\circ\big((a_1,a_2),b,c\big) := \DR^\circ\big(\PT_{a_1}a_2,b,c\big) - \P_{a_1}\DR^\circ(a_2,b,c)
	$$
	is continuous from $\CC^{\alpha_1}\times\CC^{\alpha_2}\times L^\infty\times\CC^\gamma$ to $\CC^{\alpha_1+\alpha_2+\gamma}$.   \vspace{0.15cm}

	\item Let $\beta_1,\beta_2\in(0,1)$ and $\gamma\in(-3,3)$. Then the operator
	$$
	\DR^\circ\big(a,(b_1,b_2),c\big) := \DR^\circ\big(a,\PT_{b_1}b_2,c\big) - \DR^\circ(ab_1,b_2,c)
	$$
	is continuous from $L^\infty\times\CC^{\beta_1}\times\CC^{\beta_2}\times\CC^\gamma$ to $\CC^{\beta_1+\beta_2+\gamma}$.
\end{itemize}
\end{prop}

\bigskip

\section{Correctors and commutators}
\label{AppendixContinuity}


In order to simplify the notation, we write here $\|\cdot\|_\alpha$ for $\|\cdot\|_{\mcC^\alpha}$. The proofs of the corrector estimates follow the line of reasoning of similar estimates proved in \cite{BB3}. Recall the definitions of the following operators
\begin{align*}
{\sf C}_L^<\Big(a_1,a_2,b\Big) &= \P_{L\PT_{a_1}a_2}b-a_1\P_{La_2}b,   \\
{\sf C}_L^>\Big(a,b_1,b_1\Big) &= \P_{La}\Big(\PT_{b_1}b_2\Big)-b_1\P_{La}b_2,   \\
{\sf C}_L\Big(a_1,a_2,b\Big) &= \PI\Big(L\PT_{a_1}a_2,b\Big)-a_1\PI\Big(La_2,b\Big).
\end{align*}

\medskip

\begin{thm}        \label{ThmCL}
\begin{itemize}
	\item Let $\alpha_1\in(0,1)$ and $\alpha_2,\beta\in(-3,3)$ such that $\alpha_1+\alpha_2\in(-3,3)$. If
	$$
	\alpha_2+\beta-2<0\quad\text{and}\quad\alpha_1+\alpha_2+\beta-2>0
	$$
	then the operators ${\sf C}_L^<$ and ${\sf C}_L$ have natural extensions as continuous operators from $\mcC^{\alpha_1}\times\mcC^{\alpha_2}\times\mcC^\beta$ to $\mcC^{\alpha_1+\alpha_2+\beta-2}$.
	\item Let $\beta_1\in(0,1)$ and $\alpha,\beta_2\in(-3,3)$ such that $\beta_1+\beta_2\in(-3,3)$. If
	$$
	\alpha+\beta_2-2<0\quad\text{and}\quad\alpha+\beta_1+\beta_2-2>0
	$$
	then the operator ${\sf C}_L^>$ has a natural extension as a continuous operator from $\mcC^\alpha\times\mcC^{\beta_1}\times\mcC^{\beta_2}$ to $\mcC^{\alpha+\beta_1+\beta_2-2}$.
	\item Let $\alpha_1\in(0,1)$ and $\alpha_2,\beta\in(-3,3)$ such that $\alpha_1+\alpha_2\in(-3,3)$. If
	$$
	\alpha_2+\beta-1<0\quad\text{and}\quad\alpha_1+\alpha_2+\beta-1>0
	$$
	then the operators ${\sf C}_{V_i}^<$ and ${\sf C}_{V_i}$ have natural extensions as continuous operators from $\mcC^{\alpha_1}\times\mcC^{\alpha_2}\times\mcC^\beta$ to $\mcC^{\alpha_1+\alpha_2+\beta-1}$.
	\item Let $\beta_1\in(0,1)$ and $\alpha,\beta_2\in(-3,3)$ such that $\beta_1+\beta_2\in(-3,3)$. If
	$$
	\alpha+\beta_2-1<0\quad\text{and}\quad\alpha+\beta_1+\beta_2-1>0
	$$
	then the operator ${\sf C}_{V_i}^>$ has a natural extension as a continuous operator from $\mcC^\alpha\times\mcC^{\beta_1}\times\mcC^{\beta_2}$ to $\mcC^{\alpha+\beta_1+\beta_2-1}$.
\end{itemize}
\end{thm}

\medskip

\begin{Dem}
We give here the details for the continuity estimate on $\DC_L$ and explain how to adapt the proof for $\DC_L^<,\DC_L^>,\DC_{V_i},\DC_{V_i}^<$ and $\DC_{V_i}^>$.

\ssk

We want to compute the regularity of $\DC_L(a_1,a_2,b)$ using a family $\CQ$ of $\mathsf{StGC}^r$ with $r>|\alpha_1+\alpha_2+\beta-2|$. Recall that a term $\PI(La,b)$ can be written as a linear combination of terms of the form
$$
\int_0^1\CP_t^{1\bullet}(\CQ_t^1(tL)a\cdot\CQ_t^2b)\frac{dt}{t^2},
$$
while $\PT_ba$ is a linear combination of terms of the form
$$
\int_0^1\widetilde\CQ_t^{3\bullet}\big(\widetilde\CQ_t^4a\cdot\CP_t^2b\big)\frac{dt}{t}
$$
with $\CQ^1,\CQ^2,\widetilde\CQ^4\in\mathsf{StGC}^{\frac{3}{2}}$, $\widetilde\CQ^3\in\mathsf{GC}^{\frac{3}{2}}$ and $\CP^1,\CP^2\in\mathsf{StGC}^{[0,3]}$. For the terms where $\CP^2\in\mathsf{StGC}^{[1,3]}$, we already have the correct regularity since
\begin{align*}
\int_0^1\int_0^1\CQ_u\CP_t^{1\bullet}&\left(\CQ_t^1(tL)\widetilde\CQ_s^{3\bullet}\left(\CP_s^2a_1\cdot\widetilde\CQ_s^4a_2\right)\cdot\CQ_t^2b\right)\frac{ds}{s}\frac{dt}{t^2}  \\
&\lesssim\|a_1\|_{\alpha_1}\|a_2\|_{\alpha_2}\|b\|_\beta\int_0^1\int_0^1\left(\frac{ut}{(t+u)^2}\right)^{\frac{r}{2}}\left(\frac{ts}{(s+t)^2}\right)^{\frac{3}{2}}s^{\frac{\alpha_1+\alpha_2}{2}}t^{\frac{\beta}{2}}\frac{ds}{s}\frac{dt}{t^2}  \\
&\lesssim\|a_1\|_{\alpha_1}\|a_2\|_{\alpha_2}\|b\|_\beta\ u^{\frac{\alpha_1+\alpha_2+\beta-2}{2}}
\end{align*}
using that $\alpha_1\in(0,1)$. We only consider $\CP^2\in\mathsf{StGC}^0$ for the remainder of the proof. For all $e\in\mcM$, we have
$$
\DC_L(a_1,a_2,b)(e)=\PI\left(L\PT_{a_1}a_2,b\right)(e)-a_1(e)\cdot\PI(La_2,b)(e)=\PI\left(L\PT_{a_1}a_2-a_1(e)\cdot La_2,b\right)(e),
$$
since $\PI$ is bilinear and $a_1(e)$ is a scalar. This yields that $\DC(a_1,a_2,b)(e)$ is a linear combination of terms of the form
$$
\int_0^1\int_0^1\CP_t^{1\bullet}\bigg(\CQ_t^1(tL)\widetilde\CQ_s^{3\bullet}\left(\left(\CP_s^2a_1-a_1(e)\right)\cdot\widetilde\CQ_s^4a_2\right)\cdot\CQ_t^2b\bigg)(e)\,\frac{ds}{s}\frac{dt}{t^2}
$$
using that $\dst\int_0^1L\widetilde\CQ_s^{3\bullet}\widetilde\CQ_s^4\frac{ds}{s}=L$, up to smooth terms. This gives $\big(\CQ_u\DC_L(a_1,a_2,b)\big)(e)$ as a linear combination of terms of the form
\begin{small}\begin{align*}
&\int K_{\CQ_u}(e,e') \CP_t^{1\bullet}\bigg(\CQ_t^1(tL)\widetilde\CQ_s^{3\bullet}\left(\left(\CP_s^2a_1-a_1(e')\right)\cdot\widetilde\CQ_s^4a_2\right)\cdot\CQ_t^2b\bigg)(e')\,\frac{ds}{s}\frac{dt}{t^2}\nu(de')   \\
&= \int K_{\CQ_u}(e,e')K_{\CP_t^{1\bullet}}(e',e'')\bigg(\CQ_t^1(tL)\widetilde\CQ_s^{3\bullet}\left(\left(\CP_s^2a_1-a_1(e'')\right)\cdot\widetilde\CQ_s^4a_2\right)\cdot\CQ_t^2b\bigg)(e'')\frac{ds}{s}\frac{dt}{t^2}\nu(de')\nu(de'')   
\end{align*}
\begin{align*}
&\quad +\int\int_0^u K_{\CQ_u}(e,e')K_{\CP_t^{1\bullet}}(e',e'')\Big(a_1(e'')-a_1(e')\Big)\left(\CQ_t^1(tL)a_2\cdot\CQ_t^2b\right)(e'')\frac{dt}{t^2}\nu(de')\nu(de'')   \\
&\quad +\int\int_u^1 K_{\CQ_u}(e,e')K_{\CP_t^{1\bullet}}(e',e'')\Big(a_1(e'')-a_1(e')\Big)\left(\CQ_t^1(tL)a_2\cdot\CQ_t^2b\right)(e'')\frac{dt}{t^2}\nu(de')\nu(de'')   \\
&=: A+B+C.
\end{align*}\end{small}
The term $A$ is bounded using cancellations properties. We have

\begin{align*}
|A|&= \int K_{\CQ_u\CP_t^{1\bullet}}(e,e')\bigg(\CQ_t^1(tL)\widetilde\CQ_s^{3\bullet}\left(\left(\CP_s^2a_1-a_1(e')\right)\cdot\widetilde\CQ_s^4a_2\right)\cdot\CQ_t^2b\bigg)(e') \frac{ds}{s}\frac{dt}{t^2}\nu(de')   \\
&\lesssim \|a_1\|_{\alpha_1}\|a_2\|_{\alpha_2}\|b\|_\beta\left(\int_0^u\int_0^1\left(\frac{st}{(s+t)^2}\right)^{\frac{3}{2}}(s+t)^{\frac{\alpha_1}{2}}s^{\frac{\alpha_2}{2}}t^{\frac{\beta}{2}}\frac{ds}{s}\frac{dt}{t^2}\right.\\
&\qquad \left.+\int_u^1\int_0^1\left(\frac{tu}{(t+u)^2}\right)^{\frac{r}{2}}\left(\frac{st}{(s+t)^2}\right)^{\frac{3}{2}}(s+t)^{\frac{\alpha_1}{2}}s^{\frac{\alpha_2}{2}}t^{\frac{\beta}{2}}\frac{ds}{s}\frac{dt}{t^2}\right)\\
&\lesssim \|a_1\|_{\alpha_1}\|a_2\|_{\alpha_2}\|b\|_\beta\ u^{\frac{\alpha_1+\alpha_2+\beta-2}{2}},
\end{align*}
using that $\alpha_1\in(0,1),\CP^2\in\mathsf{StGC}^0$ and $(\alpha_1+\alpha_2+\beta-2)>0$. 


For the term $B$, we have
\begin{align*}
|B|&\lesssim\|a_1\|_{\alpha_1}\|a_2\|_{\alpha_2}\|b\|_\beta\int_{e',e''}\int_0^uK_{\CQ_u}(e,e')K_{\CP_t^{1\bullet}}(e',e'')d(e',e'')^{\alpha_1}t^{\frac{\alpha_2+\beta}{2}}\frac{dt}{t^2}\nu(de')\nu(de'')\\
&\lesssim\|a_1\|_{\alpha_1}\|a_2\|_{\alpha_2}\|b\|_\beta\int_0^ut^{\frac{\alpha_1+\alpha_2+\beta-2}{2}}\frac{dt}{t}\\
&\lesssim\|a_1\|_{\alpha_1}\|a_2\|_{\alpha_2}\|b\|_\beta\ u^{\frac{\alpha_1+\alpha_2+\beta-2}{2}},
\end{align*}
using again that $\alpha_1\in(0,1)$ and $(\alpha_1+\alpha_2+\beta-2)>0$. 


Finally for $C$, we also use cancellations properties to get
\begin{align*}
|C|&\lesssim\|a_1\|_{\alpha_1}\|a_2\|_{\alpha_2}\|b\|_\beta \bigg\{\int_{e',e''}\int_u^1K_{\CQ_u}(e,e')K_{\CP_t^{1\bullet}}(e',e'')\Big|a_1(e)-a_1(e')\Big|t^{\frac{\alpha_2+\beta}{2}}\frac{dt}{t^2}\nu(de')\nu(de'')   \\
&\quad+ \int_{e',e''}\int_u^1K_{\CQ_u}(e,e')K_{\CP_t^{1\bullet}}(e',e'')\Big|a_1(e')-a_1(e'')\Big|t^{\frac{\alpha_2+\beta}{2}}\frac{dt}{t^2}\nu(de')\nu(de'')\bigg\}   \\
&\lesssim\|a_1\|_{\alpha_1}\|a_2\|_{\alpha_2}\|b\|_\beta \bigg\{\int_{e',e''}\int_u^1K_{\CQ_u}(e,e')K_{\CP_t^{1\bullet}}(e',e'')d(e,e')^{\alpha_1}t^{\frac{\alpha_2+\beta}{2}}\frac{dt}{t^2}\nu(de')\nu(de'')   \\
&\quad+ \int_{e',e''}\int_u^1K_{\CQ_u}(e,e')K_{\CP_t^{1\bullet}}(e',e'')d(e',e'')^{\alpha_1}t^{\frac{\alpha_2+\beta}{2}}\frac{dt}{t^2}\nu(de')\nu(de'')\bigg\}   
\end{align*}
\begin{align*}
&\lesssim\|a_1\|_{\alpha_1}\|a_2\|_{\alpha_2}\|b\|_\beta \bigg\{u^{\frac{\alpha_1}{2}}\int_u^1t^{\frac{\alpha_2+\beta-2}{2}}\frac{dt}{t}+\int_u^1 \left(\frac{tu}{(t+u)^2}\right)^{\frac{r}{2}}t^{\frac{\alpha_1+\alpha_2+\beta-2}{2}}\frac{dt}{t}\bigg\}  \\
&\lesssim\|a_1\|_{\alpha_1}\|a_2\|_{\alpha_2}\|b\|_\beta\ u^{\frac{\alpha_1+\alpha_2+\beta-2}{2}},
\end{align*}
using that $\alpha_1\in(0,1)$ and $(\alpha_2+\beta-2)<0$. In the end, we have
$$
\Big\|\CQ_u\DC(a_1,a_2,b)\Big\|_\infty\lesssim\|a_1\|_{\alpha_1}\|a_2\|_{\alpha_2}\|b\|_\beta\ u^{\frac{\alpha_1+\alpha_2+\beta-2}{2}}
$$
uniformly in $u\in(0,1]$, so the proof is complete for $\DC_L$. The proofs for $\DC_L^<$ and $\DC_L^>$ are then easy to obtain since $\P_{La}b$ has the same form as $\PI(La,b)$. Indeed, $\P_{La}b$ is a linear combination of
$$
\int_0^1\CQ_t^{1\bullet}\Big(\CP_t^1(tL)a\cdot\CQ_t^2b\Big)\,\frac{dt}{t^2}
$$
where $\CQ^1,\CQ^2\in\mathsf{StGC}^{\frac{3}{2}},\CP^1\in\mathsf{StGC}^{[0,3]}$ and we have $\big(\CP_t^1(tL)\big)_{0<t\le1}\in\mathsf{StGC}^2$. 

\ssk

The proofs for $\DC_{V_i},\DC_{V_i}^<$ and $\DC_{V_i}^>$ also follow from the same argument and using the Leibnitz rule as for the corrector $\DC_\partial$ use in \cite{BB3} to solve the generalised (KPZ) equation.
\end{Dem}

\medskip

\begin{thm}        \label{ThmCL1}
\begin{itemize}
    \item Let $\alpha_1\in(1,2)$ and $\alpha_2,\beta\in(-3,3)$ such that $\alpha_1+\alpha_2\in(-3,3)$. If
    $$\alpha_2+\beta-2<0\quad\text{and}\quad\alpha_1+\alpha_2+\beta-2>0$$
    then the operators $\DC_{L,(1)}^<$ and $\DC_{L,(1)}$ have natural extension as continuous operators from $\CC^{\alpha_1}\times\CC^{\alpha_2}\times\CC^\beta$ to $\CC^{\alpha_1+\alpha_2+\beta-2}$.
    \item Let $\beta_1\in(1,2)$ and $\alpha,\beta_2\in(-3,3)$ such that $\beta_1+\beta_2\in(-3,3)$. If
    $$\alpha+\beta_2-2<0\quad\text{and}\quad\alpha+\beta_1+\beta_2-2>0$$
    then the operator $\DC_{L,(1)}^>$ has a natural extension as a continuous operator from $\CC^\alpha\times\CC^{\beta_1}\times\CC^{\beta_2}$ to $\CC^{\alpha+\beta_1+\beta_2-2}$.
\end{itemize}
\end{thm}

\medskip

\begin{Dem}
For the continuity estimate of $\DC_{L,(1)}$, we also want to compute the regularity using a family $\CQ$ of $\mathsf{StGC}^r$ with $r>|\alpha_1+\alpha_2+\beta-2|$. Again a term $\PI(La,b)$ can be written as a linear combination of terms of the form
$$
\int_0^1\CP_t^{1\bullet}\big(\CQ_t^1(tL)a\cdot\CQ_t^2b\big)\frac{dt}{t^2},
$$
while $\PT_ba$ is a linear combination of terms of the form
$$
\int_0^1\widetilde\CQ_t^{3\bullet}(\widetilde\CQ_t^4a\cdot\CP_t^2b)\frac{dt}{t},
$$
with $\CQ^1,\CQ^2,\widetilde\CQ^3,\widetilde\CQ^4\in\mathsf{StGC}^{\frac{3}{2}}$ and $\CP^1,\CP^2\in\mathsf{StGC}^{[0,3]}$. For the terms where $\CP^2\in\mathsf{StGC}^{[2,3]}$, we already have the correct regularity since
\begin{align*}
\int_0^1\int_0^1\CQ_u\CP_t^{1\bullet}&\left(\CQ_t^1(tL)\widetilde\CQ_s^{3\bullet}\left(\CP_s^2a_1\cdot\widetilde\CQ_s^4a_2\right)\cdot\CQ_t^2b\right)\frac{ds}{s}\frac{dt}{t^2}\\
&\lesssim\|a_1\|_{\alpha_1}\|a_2\|_{\alpha_2}\|b\|_\beta\int_0^1\int_0^1\left(\frac{ut}{(t+u)^2}\right)^{\frac{r}{2}}\left(\frac{ts}{(s+t)^2}\right)^{\frac{3}{2}}s^{\frac{\alpha_1+\alpha_2}{2}}t^{\frac{\beta}{2}}\frac{ds}{s}\frac{dt}{t^2}\\
&\lesssim\|a_1\|_{\alpha_1}\|a_2\|_{\alpha_2}\|b\|_\beta\ u^{\frac{\alpha_1+\alpha_2+\beta-2}{2}}
\end{align*}
using that $\alpha_1\in(1,2)$ so we only consider $\CP^2\in\mathsf{StGC}^{[0,1]}$. For $\CP^2\in\mathsf{StGC}^0$, we control it using the term $a_1\PI(La_2,b)$ as in the proof of the continuity estimate of $\DC$. We are left with
$$
\int\CP_t^{1\bullet}\left(\CQ_t^1(tL)\widetilde\CQ_s^{3\bullet}\left(\left(\CP_s^2\Big(a_1-d\big(\overline{u}_0(e)\big)^{-1}\sum_{i=1}^\ell(V_ia_1)(e)\delta_i(\cdot,e)\Big)\right)\cdot\widetilde\CQ_s^4a_2\right)\cdot\CQ_t^2b\right)(e)\frac{ds}{s}\frac{dt}{t^2}
$$
with $\CP^2\in\mathsf{StGC}^1$. Then the result follows with the same proof using that $\CP_s^21=0$ since it encodes some cancellation and the first order Taylor expansion
$$
\left|a_1(e')-a_1(e)-d(\overline{u}_0)^{-1}\sum_{i=1}^\ell(V_ia_1)(e)\delta_i(e',e)\right|\lesssim d(e,e')^\alpha.
$$
We let the reader prove the continuity resuls for $\DC_{L,(1)}^<$ and $\DC_{L,(1)}^>$; they can be proved by the same argument as above.
\end{Dem}

\medskip


\begin{thm}        \label{ThmLV}
\begin{itemize}
    \item Let $\alpha\in(0,1)$ and $\beta\in(-3,3)$, be such that $\alpha+\beta<3$, and $\alpha+\beta-2\in(-3,3)$. Then the operator ${\DL}$ has a natural extension as a continuous operator from $\CC^{\alpha}\times\CC^{\beta}$ into $\CC^{\alpha+\beta-2}$.   \vspace{0.1cm}
    
    \item Let $\alpha_1,\alpha_2\in(0,1)$ and $\beta\in(-3,3)$ such that $\alpha_1+\beta<3$ and $\alpha_1+\alpha_2+\beta-2\in(-3,3)$. Then the iterated operator 
    $$
    {\DL}\big((a_1,a_2),b\big) := {\DL}\big(\P_{a_1}a_2,b\big) - \P_{a_1}{\DL}(a_2,b)
    $$ 
    has a natural extension as a continuous operator from $\CC^{\alpha_1}\times\CC^{\alpha_2}\times\CC^\beta$ into $\CC^{\alpha_1+\alpha_2+\beta-2}$.   \vspace{0.1cm}
    
    \item Let $\alpha_1,\alpha_2,\alpha_3\in(0,1)$ and $\beta\in(-3,3)$ such that $\alpha_1+\alpha_2+\beta<3$, $\alpha_2+\beta<3$ and $\alpha_1+\alpha_2+\beta-2\in(-3,3)$. Then the iterated operator 
    $$
    {\DL}\Big(\big((a_1,a_2),a_3\big), b\big) := {\DL}\big((\P_{a_1}a_2,a_3), b\big) - \P_{a_1}{\DL}\big((a_2,a_3), b\big)
    $$ 
    has a natural extension as a continuous operator from $\CC^{\alpha_1}\times\CC^{\alpha_2}\times\CC^{\alpha_3}\times\CC^\beta$ into $\CC^{\alpha_1+\alpha_2+\alpha_3+\beta-2}$.\vspace{0.1cm}

    \item Let $\alpha\in(1,2)$ and $\beta\in(-3,3)$, be such that $\alpha+\beta<3$, and $\alpha+\beta-2\in(-3,3)$. Then the operator ${\DL_{(1)}}$ has a natural extension as a continuous operator from $\CC^{\alpha}\times\CC^{\beta}$ into $\CC^{\alpha+\beta-2}$.

    \item Let $\alpha,\beta\in(-3,3)$ such that $\alpha+\beta-1\in(-3,3)$. Then the operator $\DV_i$ has a natural extension as a continuous operator from $\mcC^\alpha\times\mcC^\beta$ to $\mcC^{\alpha+\beta-1}$.\vspace{0.1cm}
    
    \item Let $\alpha_1,\alpha_2\in(0,1)$ and $\beta\in(-3,3)$ such that $\alpha_1+\beta<3$ and $\alpha_1+\alpha_2+\beta-1\in(-3,3)$. Then the iterated operator 
    $$
    \DV_i((a_1,a_2),b):=\DV_i(\P_{a_1}a_2,b)-\P_{a_1}\DV_i(a_2,b)
    $$ 
    has a natural extension as a continuous operator from $\mcC^{\alpha_1}\times\mcC^{\alpha_2}\times\mcC^\beta$ to $\mcC^{\alpha_1+\alpha_2+\beta-1}$.
\end{itemize}
\end{thm}

\medskip

\begin{Dem}
We give the proof for the continuity estimate on $\DL$ and $\DL_{(1)}$. We let the reader adapt the proof from \cite{BB3} for the iterated operators of $\DL$ since it relies on the same argument. The same holds for $\DV_i(a,b)$ and its first iteration.

\ssk

We want to compute the regularity of $\DL(a,b)=L\PT_ab-\P_aLb$ using a family $\CQ\in\mathsf{StGC}^r$ with $r>|\alpha+\beta-2|$. We write $\PT_ab$ and $\P_ab$ respectively as linear combination of
$$
\int_0^1\widetilde\CQ_s^{3\bullet}\left(\CP_s^2a\cdot\widetilde\CQ_s^4b\right)\frac{ds}{s}\quad\text{and}\quad\int_0^1\CQ_t^{1\bullet}\left(\CP_t^1a\cdot\CQ_t^2b\right)\frac{dt}{t}
$$
with $\CQ^1,\CQ^2,\widetilde\CQ^4\in\mathsf{StGC}^{\frac{3}{2}},\widetilde\CQ^3\in\mathsf{GC}^{\frac{3}{2}}$ and $\CP^1,\CP^2\in\mathsf{StGC}^{[0,3]}$. As done for $\DC$, we only have to consider $\CP^1,\CP^2\in\mathsf{StGC}^0$ since the other terms already have the right regularity using that $\alpha\in(0,1)$. We consider a term
$$
\int_0^1L\widetilde\CQ_s^{3\bullet}\left(\CP_s^2a\cdot\widetilde\CQ_s^4b\right)\frac{ds}{s}-\int_0^1\CQ_t^{1\bullet}\left(\CP_t^1a\cdot\CQ_t^2(tL)b\right)\frac{dt}{t^2}.
$$
We use that $\dst\int_0^1\CQ_t^{1\bullet}\CQ_t^2\frac{dt}{t}=\int_0^1\widetilde\CQ_s^{3\bullet}\widetilde\CQ_s^4\frac{ds}{s}=\text{Id}$ up to smooth term to get
\begin{align*}
\int_0^1\int_0^1\CQ_t^{1\bullet}&\left(\CQ_t^2(tL)\widetilde\CQ_s^{3\bullet}\left(\CP_s^2a\cdot\widetilde\CQ_s^4b\right)-\CP_t^1a\cdot\CQ_t^2(tL)\widetilde\CQ_s^{3\bullet}\widetilde\CQ_s^4b\right)\frac{dt}{t^2}\frac{ds}{s}\\
&\quad=\int_0^1\int_0^1\CQ_t^{1\bullet}\left(\CQ_t^2(tL)\widetilde\CQ_s^{3\bullet}\left(\big(\CP_s^2a-\CP_t^1a(\cdot)\big)\cdot\widetilde\CQ_s^4b\right)\right)\frac{dt}{t^2}\frac{ds}{s}
\end{align*}
where the variable of $\CP_t^1a(\cdot)$ is frozen as before, in the sense that $\CQ_t^2(tL)\widetilde\CQ_s^{3\bullet}$ does not act on it. Since $\alpha\in(0,1)$, we can use that for any $e,e'$
$$
\left|\CP_s^2a(e')-\CP_t^1a(e)\right|\le\left|\CP_s^2a(e')-a(e')\right|+|a(e')-a(e)|+\left|a(e)-\CP_t^1a(e)\right|,
$$
to get
\begin{align*}
\int_0^1\int_0^1\CQ_u\CQ_t^{1\bullet}&\left(\CQ_t^2(tL)\widetilde\CQ_s^{3\bullet}\left(\big(\CP_s^2a-\CP_t^1a(\cdot)\big)\cdot\widetilde\CQ_s^4b\right)\right)\frac{dt}{t^2}\frac{ds}{s}\\
&\lesssim\|a\|_\alpha\|b\|_\beta\int_0^1\int_0^1\left(\frac{tu}{(t+u)^2}\right)^{\frac{r}{2}}\left(\frac{st}{(s+t)^2}\right)^{\frac{3}{2}}(t+s)^\alpha s^\beta\frac{dt}{t^2}\frac{ds}{s}\\
&\lesssim\|a\|_\alpha\|b\|_\beta\ u^{\frac{\alpha+\beta}{2}}
\end{align*}
which complete the proof for $\DL(a,b)$.

\medskip

We finally prove the estimate for the refined commutator $\DL_{(1)}(a,b)$ that is given for any $e\in\CM$ by
$$
\DL_{(1)}(a,b)(e)=\big(L\PT_ab\big)(e)-\big(\P_aLb\big)(e)-\sum_{i=1}^\ell\big(\P_{d(\overline{u}_0)^{-1}V_ia)}^{(i)}b\big)(e).
$$
where 
$$
\dst\big(\P_a^{(i)}b\big)(e)=\int_{e',e''}K(e;e',e'')a(e')\left(\PT_{\delta_i(\cdot,e')}b\right)(e'')\nu(de')\nu(de''),
$$ 
with $K$ the kernel of the bilinear operator $(a,b)\mapsto\P_ab$. As in the proof of $\DC_{L,(1)}$, we are left with
\begin{align*}
&\int K_{\CQ_t^{1\bullet}}(e,e') \bigg\{\CQ_t^2(tL)\widetilde\CQ_s^{3\bullet}\big(\CP_s^2a\cdot\widetilde\CQ_s^4b\big)   \\
&\quad -\sum_{i=1}^\ell\big(\CP_t^1(d(\overline{u_0})^{-1}V_ia)\big)(e')\cdot\CQ_t^2(tL)\widetilde\CQ_s^{3\bullet}\big(\CP_s^2\delta_i(\cdot,e')\cdot\widetilde\CQ_s^4b\big)\bigg\}(e') \, \frac{dt}{t^2}\frac{ds}{s}\nu(de')   \\
&= \int K_{\CQ_t^{1\bullet}}(e,e')\left(\CQ_t^2(tL)\widetilde\CQ_s^{3\bullet} \left(\CP_s^2\Big(a-\sum_{i=1}^\ell\CP_t^1\big(d(\overline{u}_0)a\big)(e')\big)\delta_i(\cdot,e')\Big)\cdot\widetilde\CQ_s^4b\right)\right)(e') \, \frac{dt}{t^2}\frac{ds}{s}\nu(de')
\end{align*}
with $\CP^1,\CP^2\in\mathsf{StGC}^1$. The result follows with the same proof using that $\CP_s^21=0$ since it encodes some cancellation and the first order Taylor expansion for $a$.
\end{Dem}

\vfill \pagebreak

\section{Paracontrolled expansion}
\label{SectionExpansionFormula}


\begin{thm}
Let $f:\IR\to\IR$ be a $C^4$ function and let $u$ and $v$ be respectively $C^\alpha$ and $C^{4\alpha}$ functions on $[0,T]\times\IT^3$ with $\alpha\in(0,1)$. Then
\begin{align*}
f(u)v&=\P_{f'(u)v}u+\frac{1}{2}\Big\{\P_{f^{(2)}(u)v}u^2-2\P_{f^{(2)}(u)uv}u\Big\}\\
&\quad+\frac{1}{3!}\Big\{\P_{f^{(3)}(u)v}u^3-3\P_{f^{(3)}(u)uv}u^2+3\P_{f^{(3)}(u)u^2v}u\Big\}+f_v(u)^\sharp
\end{align*}
for some remainder $f_v(u)^\sharp\in\mcC^{4\alpha}$.
\end{thm}

\medskip

\begin{Dem}
We have to prove that
\begin{align*}
R&:=vf(u)-\P_{vf'(u)}u-\frac{1}{2}\Big\{\P_{vf^{(2)}(u)}u^2-2\P_{vf^{(2)}(u)u}u\Big\}\\
&\quad-\frac{1}{3!}\Big\{\P_{vf^{(3)}(u)}u^3-3\P_{vf^{(3)}(u)u}u^2+3\P_{vf^{(3)}(u)u^2}u\Big\}
\end{align*}
is a $3\alpha$-H\"older function. Using that $\P_1vf(u)=vf(u)$ up to smooth term and that $\P_ab$ is the sum of terms of the form
$$
\int_0^1\CQ_t^{1\bullet}(\CQ_t^2a\cdot\CP_t^1b)\frac{dt}{t}
$$
with $\CQ^1,\CQ^2\in\mathsf{StGC}^{\frac{3}{2}}$ and $\CP^1\in\mathsf{StGC}^{[0,3]}$, $R$ is a sum of terms of the form $\int_0^1\CQ_t^{1\bullet}(r_t)\frac{dt}{t}$ with
\begin{align*}
r_t&:=\CQ_t^2\Big(vf(u)\Big)-\CQ_t^2\Big(vf'(u)\Big)\CP_t^1(u)-\frac{1}{2}\CQ_t^2\Big(vf^{(2)}(u)\Big)\CP_t^1(u^2)+\CQ_t^2\Big(vf^{(2)}(u)u\Big)\CP_t^1(u
)\\
&\quad+\frac{1}{6}\CQ_t^2\Big(vf^{(3)}(u)\Big)\CP_t^1(u^3)+\frac{1}{2}\CQ_t^2\Big(vf^{(3)}(u)u\Big)\CP_t^1(u^2)-\frac{1}{2}\CQ_t^2\Big(vf^{(3)}(u)u^2\Big)\CP_t^1(u).
\end{align*}
We need to get a bound on $r_t$ in $L^\infty(\CM)$. We have for $e\in\CM$
\begin{align*}
r_t(e)&=\int_{\CM^2}K_{\CQ_t^2}(e,e')K_{\CP_t^1}(e,e'')\Big\{\Big(vf(u)\Big)(e')-\Big(vf'(u)\Big)(e')u(e'')-\frac{1}{2}\Big(vf^{(2)}(u)\Big)(e')u^2(e'')\\
&\quad+\Big(vf^{(2)}(u)u\Big)(e')u(e'')+\frac{1}{6}\Big(vf^{(3)}(u)\Big)(e')u^3(e'')+\frac{1}{2}\Big(vf^{(3)}(u)u\Big)(e')u^2(e'')\\
&\quad-\frac{1}{2}\Big(vf^{(3)}(u)u^2\Big)(e')u(e'')\Big\}\nu(de')\nu(de'').
\end{align*}
Using a Taylor expansion for $f$, we have
\begin{align*}
r_t(e)&=\int_{[0,1]^4}f^{(4)}\Big(u(e'')+s_4s_3s_2s_1\left(u(e')-u(e'')\right)\Big)s_3s_2s_1\left(u(e')-u(e'')\right)^4ds_4ds_3ds_2ds_1\\
&\quad+v(e')\Big(f(u(e''))+u(e'')f'(u(e''))+\frac{1}{2}u^2(e'')f^{(2)}(u(e''))+\frac{1}{3!}u^3(e'')f^{(3)}(u(e''))\Big)\\
&=(1)+(2).
\end{align*}
For the first term, we have
$$
(1)\le\|u\|_{\alpha}^4\ t^{\frac{4\alpha}{2}}
$$
and for the second term
$$
(2)\le\|u\|_{L^\infty}\|v\|_{4\alpha}\ t^{\frac{4\alpha}{2}}
$$
which allows us to conclude.
\end{Dem}

\bigskip
\bigskip

\noindent \textcolor{gray}{$\bullet$} {\sf I. Bailleul} -- Univ. Rennes, CNRS, IRMAR - UMR 6625, F-35000 Rennes, France   \\
\noindent {\it E-mail}: ismael.bailleul@univ-rennes1.fr   

\medskip

\noindent \textcolor{gray}{$\bullet$} {\sf A. Mouzard} --  Univ. Rennes, CNRS, IRMAR - UMR 6625, F-35000 Rennes, France   \\
{\it E-mail}: antoine.mouzard@univ-rennes1.fr

\end{document}